\documentclass[12pt]{amsart}
\usepackage{amsmath,amsfonts,amssymb,amscd}
\usepackage{latexsym}
\usepackage{multirow}
\usepackage[]{youngtab}
\usepackage[usenames]{color}
\usepackage{ifpdf}
\ifpdf
  \usepackage[pdftex]{graphicx}
  \usepackage{epstopdf}
\else
  \usepackage[dvips]{graphicx}
\fi
\usepackage{url}
\usepackage{pdflscape}
\usepackage[all]{xy}

\usepackage[
hyperindex=true,pagebackref=true,bookmarks=true,
colorlinks=true,linkcolor=blue,citecolor=red]{hyperref}

\allowdisplaybreaks[1]

\title [Refined composite invariants of torus knots via DAHA]
{Refined composite invariants of torus knots via DAHA}
\author[Ivan Cherednik]{Ivan Cherednik $^\dag$}
\author[Ross Elliot]{Ross Elliot $^{\ddag}$}
\date{\today}

\thanks{$^\dag$ Partially supported by NSF grant
DMS--1363138}
\thanks{$^{\ddag}$ Partially supported by a Troesh 
Family Graduate Fellowship 2014-15}

\address[I. Cherednik]{Department of Mathematics, UNC
Chapel Hill, North Carolina 27599, USA\\
chered@email.unc.edu\\}
\address[R. Elliot]{
California Institute of Technology, Pasadena, California 91125, 
USA\\
relliot@caltech.edu\\}

\begin{document}

\renewcommand{\tilde}{\widetilde}
\renewcommand{\hat}{\widehat}

\newcommand{\BR}{{\mathbb R}}
\newcommand{\BQ}{{\mathbb Q}}
\newcommand{\BC}{{\mathbb C}}
\newcommand{\BP}{{\mathbb P}}
\newcommand{\BZ}{{\mathbb Z}}
\newcommand{\BN}{{\mathbb N}}
\newcommand{\BS}{{\mathbb S}}

\newcommand{\cH}{{\mathcal H}}
\newcommand{\cA}{{\mathcal A}}
\newcommand{\cB}{{\mathcal B}}
\newcommand{\ccF}{{\mathfrak F}}
\newcommand{\cD}{{\mathcal D}}
\newcommand{\cL}{{\mathcal L}}
\newcommand{\cF}{{\mathcal F}}
\newcommand{\cP}{{\mathcal P}}
\newcommand{\cX}{{\mathcal X}}
\newcommand{\cY}{{\mathcal Y}}
\newcommand{\cS}{{\mathcal S}}
\newcommand{\cSol}{\hbox{$\mathcal Sol$}}
\newcommand{\cT}{\hbox{$\mathcal T$}}

\newcommand{\Z}{{\mathbb Z}}
\newcommand{\Q}{{\mathbb Q}}
\newcommand{\N}{{\mathbb N}}
\newcommand{\C}{{\mathbb C}}
\newcommand{\R}{{\mathbb R}}
\newcommand{\X}{{\mathbb X}}
\newcommand{\Y}{{\mathbb Y}}

\newcommand{\CH}{{\mathcal H}}
\newcommand{\CA}{{\mathcal A}}

\def\HH{\mbox{${\mathcal H}$\kern-5.2pt${\mathcal H}$}}

\newcommand{\binomial}[2]{\genfrac{(}{)}{0pt}{}{ #1 }{ #2 }}
\newcommand{\qbinomial}[2]{\genfrac{[}{]}{0pt}{}{ #1 }{ #2 }_q }
\newcommand{\qbinom}[3]{\genfrac{[}{]}{0pt}{}{ #1 }{ #2 }_{ #3 } }


\def\der{\partial}
\def\tensor{\otimes}
\def\gam{\gamma} \def\Gam{\Gamma}
\def\del{\delta} \def\Del{\Delta}
\def\kap{\kappa}
\def\lam{\lambda} \def\Lam{\Lambda}
\def\Comp{{\mathbb C}}
\def\sM{{\mathcal M}}

\newtheorem{theorem}{Theorem}[section]
\newtheorem{maintheorem}[theorem]{Main Theorem}
\newtheorem{proposition}[theorem]{Proposition}
\newtheorem{definition}[theorem]{Definition}
\newtheorem{lemma}[theorem]{Lemma}
\newtheorem{corollary}[theorem]{Corollary}
\newtheorem{notation}[theorem]{Notation}
\newtheorem{remark}[theorem]{Remark}
\newtheorem{example}[theorem]{Example}

\newtheorem{theorem }{Theorem}[section]
\newtheorem{maintheorem }[theorem]{Main Theorem}
\newtheorem{proposition }[theorem]{Proposition}
\newtheorem{definition }[theorem]{Definition}
\newtheorem{lemma }[theorem]{Lemma}
\newtheorem{corollary }[theorem]{Corollary}
\newtheorem{notation }[theorem]{Notation}
\newtheorem{remark }[theorem]{Remark}
\newtheorem{example }[theorem]{Example}

\newtheorem{ maintheorem }[theorem]{Main Theorem}
\newtheorem{ theorem}{Theorem}[section]
\newtheorem{ proposition}[theorem]{Proposition}
\newtheorem{ definition}[theorem]{Definition}
\newtheorem{ lemma}[theorem]{Lemma}
\newtheorem{ corollary}[theorem]{Corollary}
\newtheorem{ notation}[theorem]{Notation}
\newtheorem{ remark}[theorem]{Remark}
\newtheorem{ example}[theorem]{Example}

\newtheorem{thm}{Theorem}[section]
\newtheorem{prop}[thm]{Proposition}
\newtheorem{lem}[thm]{Lemma}
\newtheorem{cor}[thm]{Corollary}
\newtheorem{conj}[thm]{Conjecture}
\newtheorem{con}[thm]{Conjecture}
\newtheorem{dfn}[thm]{Definition}
\newtheorem{df}[thm]{Definition}
 \newcommand{\rem}{{\bf Comment.\ }}
 \newcommand{\rmk}{{\bf Comment.\ }}
 \newcommand{\exmp}{{\bf Example.\ }}
 \newcommand{\ex}{{\bf Example.\ }}
 \newcommand{\prob}{{\bf Problem.\ }}

\newtheorem{note}{Note} 
\renewcommand{\thenote}{}
\newtheorem*{acka}{Acknowledgments}
\newtheorem{ack}{Acknowledgments}
\renewcommand{\theack}{}
\renewcommand{\appendixname}{\bf Appendix}
\renewcommand{\proof}{{\em Proof.\ }}

\hyphenation{
ap-pen-dix as-ymp-tot-ic at-trib-uted at-trib-ut-able
Bry-li-n-sky com-mu-ta-tion de-ge-ne-rate
de-riv-a-tive dis-trib-ute equi-vari-ant ex-tra-or-di-nary  
geo-met-ric griev-ance griev-ous grad-ed ho-lo-no-my ho-mo-thetic
in-fin-ite-ly in-fin-i-tes-i-mal Ha-rish Cha-n-dra mul-ti-plic-able 
non-euclid-ean non-iso-mor-phic non-smooth par-a-digm 
par-a-bol-ic pa-rab-o-loid pa-ram-e-trize phe-nom-e-non 
post-script pseu-do-dif-fer-en-tial pseu-do-fi-nite 
qua-drat-ics quad-ra-ture Han-kel rec-tan-gle semi-def-i-nite 
set-up wide-spread Euler-ian Feb-ru-ary Gauss-ian Grothen-dieck 
Hamil-ton-ian Her-mi-t-ian her-mi-t-ian Jan-u-ary 
Japan-ese Ka-shi-wa-ra Kor-te-weg Le-gendre No-vem-ber Rie-mann-ian 
Sep-tem-ber Za-mo-lo-d-chi-kov Kni-zh-nik quan-tum Op-dam
Mac-do-nald Ca-lo-ge-ro Su-ther-land Mo-ser 
Ol-sha-net-sky  Pe-re-lo-mov in-de-pen-dent ope-ra-tors 
cy-clo-to-mic ra-tio-nal de-gen-er-a-tion 
in-ter-est-ing de-for-ma-tions de-for-ma-tion pro-ce-dure 
fol-lows ope-ra-tors  pre-serve suf-fices ap-proach 
for-mu-las con-sider its com-ple-tion cor-re-spond-ing 
au-to-mor-phism be-cause pro-por-tional fi-nal-ly let-ting 
equi-v-a-lence ge-n-er-al-ized Mac-do-nald iden-ti-ties 
cor-re-s-pond sub-dia-grams par-ti-tion na-t-u-ral-ly 
or-dered stan-dard de-for-ma-tion ar-gu-ment com-bined 
sphe-r-i-cal rep-re-sen-ta-tions tri-go-no-me-t-ric
ge-n-er-al-ly speak-ing pri-m-it-ive ir-re-du-cible 
sum-ma-tion  rep-re-sen-ta-tives pro-por-ti-o-na-li-ty
ultra-sphe-ri-cal Ro-gers}

\def\ffor{\quad\hbox{ for }\quad}
\def\wwhen{\quad\hbox{ when }\quad}
\def\wwhere{\quad\hbox{ where }\quad}
\def\aand{\quad\hbox{ and }\quad}
\def\for{\  \hbox{ for } \ }
\def\iif{ \ \hbox{ if } \ }
\def\when{ \ \hbox{ when } \ }
\def\where{\  \hbox{ where } \ }
\def\and{\  \hbox{ and } \ }
\def\and{\  \hbox{ and } \ }
\def\oor{\  \hbox{ or } \ }
\def\proof{{\em Proof. \  }}

\def\equal{\stackrel{\,\mathbf{def}}{= \kern-3pt =}}

\def\la{\lambda}
\def\La{\Lambda}
\def\om{\omega}
\def\Om{\Omega}
\def\Th{\Theta}
\def\th{\theta}
\def\al{\alpha}
\def\be{\beta}
\def\ga{\gamma}
\def\ep{\epsilon}
\def\up{\upsilon}
\def\Up{\Upsilon}
\def\de{\delta}
\def\De{\Delta}
\def\ka{\kappa}
\def\kapp{\hbox{\bf \ae}}
\def\si{\sigma}
\def\Si{\Sigma}
\def\Ga{\Gamma}
\def\ze{\zeta}
\def\io{\iota}
\def\bio{b^\iota}
\def\aio{a^\iota}
\def\twio{\tilde{w}^\iota}
\def\hwio{\hat{w}^\iota}
\def\gio{\g^\iota}
\def\Bio{B^\iota}

\def\del{\delta}
\def\pa{\partial}
\def\vp{\varphi}
\def\ve{\varepsilon}
\def\inf{\infty}

\def\vph{\varphi}
\def\vps{\varpsi}
\def\vPh{\varPhi}
\def\vep{\varepsilon}
\def\vpi{{\varpi}}
\def\vth{{\vartheta}}
\def\vsi{{\varsigma}}
\def\vrh{{\varrho}}

\def\bph{\bar{\phi}}
\def\bsi{\bar{\si}}
\def\bvp{\bar{\varphi}}

\newcommand{\bS}{{\mathbf S}}
\newcommand{\bH}{{\mathbf H}}
\newcommand{\bF}{{\mathbf F}}
\newcommand{\bE}{{\mathbf E}}

\def\tal{\tilde{\alpha}}
\def\tbe{\tilde{\beta}}
\def\tde{\tilde{\delta}}
\def\tpi{\tilde{\pi}}
\def\txi{\tilde{\xi}}
\def\tPi{\tilde{\Pi}}
\def\tPhi{\tilde{\Phi}}
\def\tV{\tilde{V}}
\def\tJ{\tilde{J}}
\def\tla{\tilde{\lambda}}
\def\tga{\tilde{\gamma}}
\def\tGa{\tilde{\Gamma}}
\def\tvs{\tilde{{\varsigma}}}
\def\tu{\tilde{u}}
\def\tU{\tilde{U}}
\def\tw{\widetilde w}
\def\tW{\widetilde W}
\def\tB{\tilde B}
\def\tv{\tilde v}
\def\tV{\tilde V}
\def\tz{\tilde z}
\def\tb{\tilde b}
\def\ta{\tilde a}
\def\tih{\tilde h}
\def\trh{\tilde {\rho}}
\def\tx{\tilde x}
\def\tf{\tilde f}
\def\tg{\tilde g}
\def\tG{\tilde G}
\def\tk{\tilde k}
\def\tl{\tilde l}
\def\tL{\tilde L}
\def\tD{\tilde D}
\def\tR{\tilde R}
\def\tP{\tilde P}
\def\tH{\tilde H}
\def\tp{\tilde p}

\def\hH{\hat{H}}
\def\hh{\hat{h}}
\def\hR{\hat{R}}
\def\hY{\hat{Y}}
\def\hX{\hat{X}}
\def\hP{\hat{P}}
\def\hT{\hat{T}}
\def\hV{\hat{V}}
\def\hG{\hat{G}}
\def\hF{\hat{F}}
\def\hw{\widehat{w}}
\def\hW{\widehat{W}}
\def\hu{\hat{u}}
\def\hs{\hat{s}}
\def\hv{\hat{v}}
\def\hb{\hat{b}}
\def\hB{\widehat{B}}
\def\hze{\hat{\zeta}}
\def\hsi{\hat{\sigma}}
\def\hrh{\hat{\rho}}
\def\hth{\hat{\theta}}
\def\hy{\hat{y}}
\def\hx{\hat{x}}
\def\hz{\hat{z}}
\def\hg{\hat{g}}
\def\he{\hat{e}}
\def\hE{\widehat{E}}

\def\B{\mathbf{B}}
\def\I{\mathbf{I}}
\def\P{\mathbf{P}}
\def\G{\mathbf{G}}
\def\S{\mathbf{S}}
\def\F{\mathbf{F}}
\def\one{\mathbf{1}}
\def\Sn{\mathbf{S}_n}
\def\0{\mathbf{0}}
\def\H{\mathbf{H}}
\def\V{\mathbf{V}}

\def\f{\mathcal{F}}
\def\çF{\mathcal{F}}
\def\o{\mathcal{O}}
\def\t{\mathcal{T}}
\def\r{\mathcal{R}}
\def\l{\mathcal{L}}
\def\m{\mathcal{M}}
\def\k{\mathcal{K}}
\def\n{\mathcal{N}}
\def\d{\mathcal{D}}
\def\p{\mathcal{P}}
\def\cP{\mathcal{P}}
\def\a{\mathcal{A}}
\def\h{\mathcal{H}}
\def\c{\mathcal{C}}
\def\y{\mathcal{Y}}
\def\e{\mathcal{E}}
\def\v{\mathcal{V}}
\def\z{\mathcal{Z}}
\def\x{\mathcal{X}}
\def\s{\mathcal{S}}
\def\g{\mathcal{G}}
\def\u{\mathcal{U}}
\def\w{\mathcal{W}}
\def\i{\mathcal{I}}
\def\j{\mathcal{J}}
\def\b{\mathcal{B}}

\def\lan{\langle}
\def\llb{(\!(}
\def\ran{\rangle}
\def\rrb{)\!)}
 \def\dim{{\hbox{\rm dim}}_{\mathbb C}\,}
\def\lng{\hbox{\rm{\tiny lng}}}
\def\sht{\hbox{\rm{\tiny sht}}}
\def\sph{\hbox{\rm{\tiny sph}}}
\def\inv{\hbox{\rm{\tiny inv}}}

\def\br#1{\langle #1 \rangle}

\def\rank{\hbox{rank}}
\def\gl{\mathfrak{gl}_N}

\newcommand{\Aut}{\operatorname{Aut}}
\newcommand{\Hom}{\operatorname{Hom}}
\newcommand{\End}{\operatorname{End}}
\newcommand{\Ind}{\operatorname{Ind}}
\newcommand{\ad}{\operatorname{ad}}
\newcommand{\pr}{\operatorname{pr}}
\newcommand{\aweyl}{\tilde{\mathbb S}_n}
\newcommand{\hec}{{\mathcal H}^t_n}
\newcommand{\Func}{{\mathcal F}({\mathbb C}^n,{\mathcal H}^t_n)}
\newcommand{\tr}{\operatorname{tr}}
\newcommand{\Out}{\operatorname{Out}}
\newcommand{\Rad}{\operatorname{Rad}}
\newcommand{\Spec}{\operatorname{Spec}}
\newcommand{\id}{\operatorname{id}}
\newcommand{\Int}{\operatorname{Int}}
\newcommand{\ct} {\operatorname{ct}}

\newcommand{\rat}{{\mathbb Q}}
\newcommand{\real}{{\mathbb R}}
\newcommand{\cplx}{{\mathbb C}}
\newcommand{\zint}{{\mathbb Z}}

\newcommand{\sq}{\phantom{1}\hfill$\qed$}
\newcommand{\Rea}{\Re}
\newcommand{\Ima}{\Im}

\newcommand{\st}{\bowtie}
\newcommand{\modd}{\mbox{\,mod\,}}
\newcommand{\lr}{\langle}
\newcommand{\rr}{\rangle}
\newcommand{\eps}{\varepsilon}
\newcommand{\phk}{\phi^{(k)}}
\newcommand{\psk}{\psi^{(k)}}
\newcommand{\Res}{\mbox{Res}\;}
\newcommand{\sgn}{\mbox{sgn}}
\newcommand{\mn} {\left\{ \begin{array}{c}m\\
n\end{array}\right\}}

\def\sX{\mathscr{X}}
\def\sH{\mathscr{H}}
\def\sY{\mathscr{Y}}
\def\TT{\mathfrak{T}}
\def\JJ{\mathfrak{J}}
\def\HH{\mathfrak{H}}
\def\FF{\mathfrak{F}}
\def\GG{\mathfrak{G}}
\def\CC{\mathfrak{C}}
\def\LL{\mathfrak{L}}

\def\BB{\mathfrak{B}}
\def\AA{\mathfrak{A}}
\def\ZZ{\mathfrak{Z}}
\def\HH{\hbox{${\mathcal H}$\kern-5.2pt${\mathcal H}$}}
\def\HHH{\hbox{${\mathbb H}$\kern-4.2pt${\mathbb H}$}}
\def\tHH{\widetilde{\HH\ }}

\font\smm=msbm10 at 12pt 
\def\symbol#1{\hbox{\smm #1}}
\def\lsmash{{\symbol n}}
\def\rsmash{{\symbol o}}
\def\#{\sharp}

\font\tenbf=cmbx10
\font\tenrm=cmr10
\font\tenit=cmti10
\font\ninebf=cmbx9
\font\ninerm=cmr9
\font\nineit=cmti9
\font\eightbf=cmbx8
\font\eightrm=cmr8
\font\eightit=cmti8
\font\sevenrm=cmr7
\font\sevenbf=cmbx7


\par
{\centering 
Dedicated to  Vadim Schechtman \\
on the occasion of his 60th birthday
\medskip
\par}

\begin{abstract}
We define composite DAHA-superpolynomials of torus
knots, depending on pairs of Young diagrams and
generalizing the composite HOMFLY-PT polynomials
in the theory of the skein of the annulus. We provide
various examples. Our superpolynomials extend the
DAHA-Jones (refined) polynomials and satisfy all 
standard symmetries of the DAHA-superpolynomials of 
torus knots. The latter are conjecturally related to the HOMFLY-PT 
homology; such a connection is a challenge in the theory of the
annulus. At the end, we construct two DAHA-hyperpolynomials
extending the DAHA-Jones polynomials of type $E$ 
and closely related to the exceptional Deligne-Gross
series of root systems; this theme is of experimental nature. 
\end{abstract}

\def\sht{\raisebox{0.4ex}{\hbox{\rm{\tiny sht}}}}
 \def\bysame{{\bf --- }}
 \def\={{\bf --}}
 \def\rr{{\mathsf r}}
 \def\ss{{\mathsf s}}
 \def\mm{{\mathsf m}}
 \def\pp{{\mathsf p}}
 \def\ll{{\mathsf l}}
 \def\aa{{\mathsf a}}
 \def\bb{{\mathsf b}}
 \def\NS{\hbox{\tiny\sf ns}}
 \def\ssum{\hbox{\small$\sum$}}
\newcommand{\comment}[1]{}
\renewcommand{\tilde}{\widetilde}
\renewcommand{\hat}{\widehat}
\renewcommand{\V}{\mathbb{V}}
\renewcommand{\F}{\mathbb{F}}
\newcommand{\dagx}{\hbox{\tiny\mathversion{bold}$\dag$}}
\newcommand{\ddagx}{\hbox{\tiny\mathversion{bold}$\ddag$}}
\newtheorem{conjecture}[theorem]{Conjecture}
\newcommand*\toeq{
\raisebox{-0.15 em}{\,\ensuremath{
\xrightarrow{\raisebox{-0.3 em}{\ensuremath{\sim}}}}\,}
}
\newcommand{\unknot}{\hbox{\tiny\!\raisebox{0.2 em}{$\bigcirc$}}}

\maketitle
\smallskip

{\tiny
\centerline{{\bf MSC} (2010): 17B22, 17B45, 20C08,
20F36, 33D52, 57M27}
}
\smallskip

\renewcommand{\baselinestretch}{1.2}.
{
\tableofcontents 
} 
\renewcommand{\baselinestretch}{1.2}.
\vfill\eject 

\renewcommand{\natural}{\wr}

\setcounter{section}{-1} \setcounter{equation}{0} 

\section{\sc Introduction}

We introduce and study the 
{\em composite DAHA-superpolynomials\,} for torus knots
and arbitrary {\em composite weights\,} \cite{K}, 
i.e. pairs of Young diagrams.  They depend on $\,a,q,t\,$ and 
unify the corresponding $n$\=series of (refined) 
DAHA-Jones $q,t$\= polynomials of type $A_n$; all 
symmetries of superpolynomials from \cite{CJ,CJJ} hold for them.
When $t=q$ and $a\mapsto -a$,  we establish their relation to the
{\em composite HOMFLY-PT polynomials\,}, studied in 
\cite{HM,MM,AM,GJKS}.

The topological composite theory is based 
on the {\em full HOMFLY-PT skein\,} of the  annulus, 
which is an algebra generated by link diagrams drawn there.
The adjoint representation is the simplest composite weight,
which connects our results with two examples of 
adjoint DAHA-superpolynomials for the Deligne-Gross exceptional
series of root systems considered at the end of the paper.

\smallskip
{\sf Topological origins.}
In the full HOMFLY-PT skein, the orientations of the components of 
the links can be simultaneously clockwise and counterclockwise 
around the annulus, which eventually results in pairs of Young 
diagrams.
It is isomorphic to the tensor square
of the ring of symmetric functions. 
The (non-full) skein has all orientations in the same 
direction, which is insufficient for the composite theory.   
The diagonalization of the {\em meridian maps\,}
in the full skein of the annulus provides a natural
and systematic way to define the composite HOMFLY-PT polynomials
for any knots and colors. 

The role of the annulus can be clearly seen in the theory 
of {\em satellite links}, which is of fundamental
value in low-dimensional topology (including our paper).  
Given a knot $K\subset S^3$ and a Young diagram, such a link is 
generally constructed from both a diagram   
$D(K)$ of $K$, called a {\em companion\,}, and a 
link diagram $Q$ in the annulus, called a {\em pattern}.
The annulus inevitably emerges here due to the
framing of $K$, an important ingredient of this construction
(which influences the output). 
\smallskip

{\sf Superpolynomials.}
The uncolored DAHA-superpolynomials of torus knots in $S^3$ are 
conjectured to coincide with the Poincar\'e polynomials 
for the reduced {\em HOMFLY-PT homology\,} or, equivalently, 
stable reduced Khovanov-Rozansky polynomials.  See e.g. 
\cite{DGR, Kh, KhR1, KhR2, Rou, Web} 
for the corresponding knot homology theories and categorification.
This is expected to hold for any rectangular Young diagrams,
though adding colors to HOMFLY-PT homology is a theoretical
and practical challenge. 
Rectangular diagrams are natural
here, since the DAHA-superpolynomials are conjecturally positive
for such diagrams and arbitrary algebraic knots. 

We note 
that the DAHA-superpolynomials
were recently defined for iterated torus knots \cite{ChD},
which includes all algebraic knots (links are in progress). 
This is a natural setting for the composite 
DAHA-superpolynomials,
but we focus here only on torus knots. 
 
\smallskip
The theory of DAHA-Jones polynomials is 
uniform for any root systems and arbitrary weights; accordingly, 
the DAHA-superpolynomials are defined for any Young diagrams. 
They are studied reasonably well by now; at least,
all conjectures about them from \cite{CJ} are verified, 
but the positivity. This is generally beyond what topology provides,
especially upon adding arbitrary colors to the theory.

\smallskip
The key open question in the composite direction we present 
concerns the relation of our composite 
DAHA-superpolynomials to HOMFLY-PT {\em homology\,} in the 
case of annulus. A theory in the annulus is in progress, and it 
seems capable of practically producing invariants for simple 
knots and colors; see \cite{QR}. However, we hesitate to conjecture 
any explicit connection because of the absence of such examples
so far. Also, the composite DAHA-superpolynomials lack 
positivity, as do those for the non-rectangular 
diagrams and non-algebraic knots. It is not clear 
how to address this challenge, though we provide some approach
of experimental nature in \cite{ChD}.

\smallskip

{\sf Exceptional series.}
We conclude this paper with
hypothetical  adjoint (quasi-minuscule)
\emph{DAHA-hyperpolynomials} 
for the torus knots $T^{3,2},T^{4,3}$ for the 
{\em exceptional ``magic" series\,}: 
$$\{e\subset A_1 \subset A_2 \subset G_2 \subset 
D_4 \subset F_4 \subset E_6 \subset E_7 \subset E_8\}$$
from \cite{DG}. This is for 
the maximal short root $\vartheta$, which is the 
highest weight of the adjoint representation.  Thus, for 
the root systems of type $A$, we make contact with the  
composite DAHA-superpolynomials.

The root systems $G_2, F_4$ are beyond our reach so far
and we managed to find such hyperpolynomials only for 
simple torus knots (though $T(4,3)$ is not too simple). 
Nevertheless, we believe that even such examples demonstrate that 
the final theory of DAHA-hyperpolynomials will eventually 
incorporate all types of root systems (not only classical).

The hyperpolynomials we found based on the functoriality from 
\cite{DG} are non-positive but have rich symmetries. We note that 
there are (quite a few) other series where the existence of the
superpolynomials can be expected, not only for those 
of Deligne-Gross type.  For instance,  we found
(joint with Sergei Gukov) the minuscule superpolynomials for  
$\{E_6,A_6,D_5\}$. Here, as for the exceptional series above, 
only small torus  knots and the simplest weights can be managed.

Generally, deg${}_a$+1 root systems are needed to determine the 
corresponding polynomial uniquely. They provide its 
evaluations at the corresponding values of $a$, which was 
used in \cite{CJ} to define superpolynomials and
hyperpolynomials for  $ABCD$ (infinite families). This is not 
the case with $E$. However, a very small number of evaluations 
appeared sufficient in the examples we managed. 
For instance, only $E_8$ and $E_7$ are needed to determine the 
exceptional DAHA-hyperpolynomial of $T^{3,2}$ 
(assuming that it satisfies some natural properties).
There is no general understanding at the moment of how to
proceed for arbitrary torus knots and weights for
exceptional root systems.  
\smallskip

{\sf The structure of the paper.}
In Section \ref{Koike}, the composite weights $[\la,\mu]$ 
(pairs of Young diagrams) and the corresponding 
representations are defined, following \cite{K}. Then we provide 
the definition of composite 
HOMFLY-PT polynomials 
$\mathcal{H}_{[\la,\mu]}(K)$ for any knot $K$ from \cite{HM}, 
via the full HOMFLY-PT skein algebra $\mathcal{C}$ of link diagrams 
in the annulus.  Finally, Proposition \ref{PropRJ}, a generalization 
of the Rosso-Jones formula, gives effective means of producing 
$\mathcal{H}_{[\la,\mu]}(T^{\rr,\ss})$ for 
$ T^{\rr,\ss}$ via. It essentially coincides with 
formula (C.6) from \cite{GJKS};
we give its proof. 

\smallskip
In Section \ref{DAHA}, we recall the main definitions and 
results from the DAHA theory used to introduce the 
DAHA-Jones (also called refined) polynomials and 
DAHA-superpolynomials from \cite{CJ,CJJ}.  
Then, we offer the main body of results of this paper. 
Theorem \ref{STABILIZ} is the existence 
(stabilization) of composite DAHA-superpolynomials and their 
evaluations  at $q=1$. Theorem \ref{SUPRDUAL} is the composite 
super-duality, 
which is proved using a reduction to the
DAHA-Jones polynomials, closely related to
the color exchange from Theorem \ref{COLOREXCH}.
The connection to the composite HOMFLY-PT polynomials
is Theorem \ref{CONNEC}.

\smallskip
Section \ref{ExC} is devoted to various examples of
composite DAHA superpolynomials and discussion of their
symmetries from the previous section. Our 
examples confirm the stabilization, connection, super-duality and 
evaluation theorems for a selection of seven composite 
partitions and simple torus knots. Section \ref{DG} is
devoted to the examples of hyperpolynomials
for the ``magic" exceptional series from \cite{DG} (the
bottom line of the triangle considered there). 
\smallskip

{\sf The key construction.}
We begin with the definition of (reduced, tilde-normalized)
\emph{DAHA-Jones polynomials}
$\widetilde{J\!D}^{R}_{\rr,\ss}(b\,;\, q,t),$ 
associated to any torus knot $T^{\rr,\ss}$, root system $R$, and 
(dominant) weight $b\in P_+$ for $R$.  This is unchanged vs.
\cite{CJ,CJJ}. We mention that they conjecturally coincide 
with the corresponding Quantum Group invariants for torus knots 
upon $t=q$ (for both $t$, in the non-simply-laced case).  This
was checked for $A_n$ for any Young diagrams in \cite{CJ} and 
in various other cases, including the formulas conjectured
there for $E_6$ (by R.~E.).

When $R$ is of type $A_n$, the DAHA-Jones polynomials 
are uniform  with respect to $n$; see \cite{CJ,GoN}. 
Namely, the corresponding superpolynomials are defined as 
follows: 
\begin{equation*}
H\!D_{\rr,\ss}(\la\,;\, q,t,a\mapsto -t^{n+1}) = 
\widetilde{J\!D}^{A_n}_{\rr,\ss}(\la\,;\, q,t),
\end{equation*}
where the Young diagram $\la$  is interpreted naturally as an 
$A_n$\=weight for any sufficiently large $n$.  This definition
is generalized in the present paper to the case
of the pairs $[\la,\mu]$ of Young diagrams, placed at the 
ends of the corresponding Dynkin graph for $A_n$.

The {\em uncolored case\,} corresponds to the
adjoint representation:
\begin{equation*}
H\!D_{\rr,\ss}([\om_1,\om_1]\,;\, q,t,a\mapsto -t^{n+1}) = 
\widetilde{J\!D}^{A_n}_{\rr,\ss}(\om_1+\om_n\,;\, q,t).
\end{equation*}
The stabilization is a more subtle issue in the composite
case. We prove that all symmetries
from \cite{CJ,CJJ} of the resulting 
{\em composite DAHA-superpolynomials\,} hold.
The key result of this paper is the
coincidence of $H\!D_{\rr,\ss}([\la,\mu]\,;\, q,q,-a)$
with the HOMFLY-PT polynomials defined in \cite{HM} for any 
composite diagrams $[\la,\mu]$ via the skein theory of link 
diagrams in the annulus $S^1\times I$. 

\setcounter{equation}{0} 
\section{\sc HOMFLY-PT polynomials}\label{HOMFLY}
\subsection{\bf Composite representations}\label{Koike}
An irreducible (finite-dimensional) representation $V$ 
of $\mathfrak{sl}_N(\mathbb{C})$ is uniquely specified by its 
{\em highest weight\,}: 
$$b = \sum_{i=1}^{N-1}b_i\omega_i \in P_+ \equal\,
\bigoplus_{i=1}^{N-1}\, 
\mathbb{Z}_{+}\,\omega_i,\
\Z_+=\Z_{\geq 0},$$
where $\{\omega_i\}$ are the fundamental dominant weights for 
$A_{N-1}$.  

Equivalently, we can represent $b$ (and $V$) 
by a {\em partition\,} or its 
corresponding {\em Young diagram\,} 
$\la=\la_1\ge\la_2\ldots\ \la_{N-1}\ge\la_N=0$ with at 
most $N-1$ nonempty rows and \,$k${\tiny\,th\,} row of length 
$\lambda_k \equal b_k+\cdots +b_{N-1}$.  The 
highest weight $b$ is recovered from $\lambda$ by 
taking $b_i = \lambda_i - \lambda_{i+1}$; i.e. $b_i$ is the
number of columns of $\la$ of height $i$.

The dual representation $V^\ast$ is specified by the highest weight 
$b^\ast \equal \iota(b)$, where $\iota:\omega_i\mapsto\omega_{N-i}$.  
Alternatively, the Young diagram $\lambda^\ast$ has rows of length 
$\lambda_k^\ast = \lambda_1 - \lambda_{N+1-k}$ (this operation
depends on $N$).

A weight $b\in P_+$ for $\mathfrak{sl}_N(\mathbb{C})$ can be 
interpreted for $\mathfrak{sl}_M(\mathbb{C})$ by setting $b_i = 0$ 
for $i\geq \min\{M,N\}$.  Accordingly, we can interpret the 
corresponding Young diagram $\lambda$ as a dominant weight for
$\mathfrak{sl}_M(\mathbb{C})$ by removing any 
columns of height $\geq M$.  It is precisely this sort of 
``packaging" of representations for all ranks that leads to the 
HOMFLY-PT polynomial and its generalizations. 

One can generalize this procedure to any number of Young
diagrams, ``placing" them in the Dynkin diagram of 
type $A_{N-1}$ with breaks in between.
The {\em composite representations\,} are labeled by pairs
of partitions (or Young diagrams) ``placed" at the ends of
the Dynkin diagram. Namely, for Young diagrams $\la$ and $\mu$
with $\ell(\la)$ and $\ell(\mu)$ rows, 
$N\ge \ell(\la)+\ell(\mu)$ (always assumed),
and $P_+$ of type $A_{N-1}$, let 
\begin{align}\label{lamuN}
[\lambda,\mu]_N=b^\ast+c\in P_+=P_+^{A_{N-1}} \for b,c
\hbox{\, associated with \,} \la,\mu.
\end{align}
We call the pair $[\la,\mu]$ a {\em composite
diagram/partition\,} and will 
constantly identify dominant weights 
$[\lambda,\mu]_N$ and the corresponding Young diagrams
(with no greater than $N-1$ rows).


\subsubsection{\sf Schur functions}
In what follows, we will require some basic facts about Schur 
functions and their generalization to composite representations 
in \cite{K}.

Let $\Lambda_n \equal \mathbb{Z}[x_1,\ldots,x_n]^{S_n}$ denote the 
\emph{ring of symmetric functions in $n$\=variables}, where 
the action of $S_n$ is by permuting the variables.  For any 
$m\geq n$, the map which sends $x_i\mapsto 0$ for $i>n$, and 
$x_i\mapsto x_i$ otherwise, is the restriction homomorphism 
$\Lambda_m\rightarrow \Lambda_n$.  Then the \emph{ring of 
symmetric functions} is
\begin{equation*}
\Lambda_x \equal \displaystyle\varprojlim_n\Lambda_n,
\end{equation*}
where the projective limit is taken with respect to the 
restriction homomorphisms.

If $\lambda$ is a partition with length at most $n$, one 
can define the corresponding \emph{Schur function} 
$s_\lambda(x_1,\ldots,x_n)\in \Lambda_n$.  The set of Schur 
functions for all such partitions is a $\mathbb{Z}$\=basis for
$\Lambda_n$.  We may naturally interpret a given 
$s_\lambda(x_1,\ldots,x_n)$ as having infinitely-many variables, 
for which we write $s_\lambda(\vec{x})\in\Lambda_x$.  The set of 
all $s_\lambda(\vec{x})$ is a $\mathbb{Z}$\=basis for $\Lambda_x$.

The Schur functions satisfy many interesting properties.  For 
our purposes, we will interpret $s_\lambda(\vec{x})\in\Lambda_x$ 
as a character for the irreducible polynomial representation 
$V_\lambda$.  Consequently, the \emph{Littlewood-Richardson rule},
that is
\begin{equation}\label{LR}
s_\lambda(\vec{x})s_\mu(\vec{x}) = \displaystyle\sum_\nu 
N_{\lambda,\mu}^\nu s_\nu(\vec{x}),
\end{equation}
shows that the multiplicity of an irreducible summand 
$V_\nu$ in the tensor product decomposition of 
$V_\lambda\otimes V_\mu$ is equal to the 
\emph{Littlewood-Richardson coefficient} $N_{\lambda,\mu}^\nu$.

\subsubsection{\sf The composite case}
In \cite{K}, the author introduces 
$s_{[\lambda,\mu]}(\vec{x},\vec{y})\in\Lambda_x\otimes\Lambda_y$,
which generalize the Schur functions and provide characters for
irreducible representations $V_{[\lambda,\mu]}$ corresponding to 
composite partitions. 
Their natural projection onto the character 
ring for $\mathfrak{sl}_N$ is the (ordinary) Schur function 
$s_{[\lambda,\mu]_N}(x_1,\ldots,x_{N-1})\in\Lambda_{N-1}$.
Recall that we always assume that 
$N\ge \ell(\la)+\ell(\mu)$ for the length $\ell(\la)$ of $\la$; 
see (\ref{lamuN}). 

The following formulas, proved in \cite{K}, will be used 
as definitions in our paper:
\begin{align}\label{KLR}
&s_{[\lambda,\mu]}(\vec{x},\vec{y}) \,\equal\,
\displaystyle\sum_{\tau,\nu,\xi} (-1)^{|\tau|}
N_{\nu,\tau}^\lambda N_{\tau,\xi}^\mu s_\nu(\vec{x})s_\xi(\vec{y}),
\\
\label{LRK}
&\hbox{where\, } s_\eta(\vec{x})s_\delta(\vec{y}) = 
\displaystyle\sum_{\alpha,\beta,\delta} 
N_{\beta,\alpha}^\eta N_{\gamma,\alpha}^\delta 
s_{[\beta,\gamma]}(\vec{x},\vec{y});
\end{align}
the sums here are over arbitrary triples of 
Young diagrams.

\subsection{\bf Skein theory in the annulus}
\subsubsection{\sf Composite HOMFLY-PT polynomials}
The colored HOMFLY-PT polynomial for a knot $K$ and a 
partition $\lambda$ is the integer Laurent polynomial 
$\mathcal{H}_\lambda(K;q,a)\in\mathbb{Z}[q^{\pm 1},a^{\pm 1}]$ 
satisfying $\mathcal{H}_\lambda(K;q;q^N) = 
\j^{\mathfrak{sl}_N}_{\lambda}(K;q)$ to the corresponding 
Jones polynomial for $\mathfrak{sl}_N$ and partition
(dominant weight) $\la$. The latter is also called
the Quantum Group knot invariant or WRT invariant.

The {\em composite HOMFLY-PT polynomial\,} 
for $[\la,\mu]$ is defined similarly via the specializations
$\mathcal{H}_{[\lambda,\mu]}(K;q,q^N) = 
\j^{\mathfrak{sl}_N}_{[\lambda,\mu]_N}(K;q)$ for all sufficiently 
large $N$.  In particular, $\mathcal{H}_{[\varnothing,\mu]}(K)= 
\mathcal{H}_\mu(K)$. Recall that the composite diagram
$[\lambda,\mu]_N$ is from (\ref{lamuN}).

The HOMFLY-PT polynomial has two normalizations.  For connection 
with DAHA, as in Theorem \ref{CONNEC}, we will be interested in the 
\emph{normalized} polynomial $\mathcal{H}$. However, for many of our 
intermediate calculations, we will also need the \emph{unnormalized} 
HOMFLY-PT polynomial $\bar{\mathcal{H}}$.  These are generally
defined and related by:
\begin{equation}\label{norm}
\bar{\mathcal{H}}(K) = \bar{\mathcal{H}}(U)\mathcal{H}(K),\ 
\bar{\mathcal{H}}(U)
=\text{dim}_{q,a}(V),
\end{equation}
where $K$ is any knot, $U$ is the unknot, and $\text{dim}_{q,a}$ is 
defined in Section \ref{qdim} for $V=V_{[\lambda,\mu]}$.
Observe that with this definition, 
$\mathcal{H}(U)=1$.  In the specializations described earlier in 
this section, the normalized (resp. unnormalized) HOMFLY-PT 
polynomials coincide with the reduced (resp. unreduced) Quantum 
Group knot invariants.
\smallskip

We will briefly recall the approach to composite 
HOMFLY-PT polynomials from \cite{HM}.  The \emph{full HOMFLY-PT 
skein algebra} $\mathcal{C}$ is a commutative algebra 
over the coefficient ring 
$\Upsilon = \mathbb{Z}[v^{\pm 1},s^{\pm 1}]
(\{s^k - s^{-k}\}_{k\geq 1})^{-1}$.
It consists in $\Upsilon$\=linear combinations of oriented 
link diagrams in 
$S^1\times I$. 

The {\em product\,} of two diagrams in $\mathcal{C}$ is the diagram 
obtained by identifying the outer circle of one annulus with the 
inner circle of the other; the identity with respect to this 
product is the empty diagram (with coefficient $1$). 

The relations in $\mathcal{C}$ are the (framed) HOMFLY-PT skein 
relation 
\[\left\langle
\def\objectstyle{\scriptstyle}
\def\labelstyle{\scriptstyle}
\vcenter{\xymatrix @=2pc {
&\\
\ar[ur]&\ar[ul]|\hole\\}}\right\rangle
 - 
\left\langle\def\objectstyle{\scriptstyle}
\def\labelstyle{\scriptstyle}
\vcenter{\xymatrix @=2pc{
&\\
\ar[ur]|\hole&\ar[ul]
}}\right\rangle
= (s - s^{-1})
\left\langle\def\objectstyle{\scriptstyle}
\def\labelstyle{\scriptstyle}
\vcenter{\xymatrix @=2pc{
&\\
\ar@/_/[u]&\ar@/^/[u]
}}\right\rangle,\]
together with the relation that accompanies a type-I Reidemeister 
move on a positively (resp. negatively) oriented loop with 
multiplication by a factor of $v^{-1}$ (resp. $v$).  As a 
consequence, observe that 
\[\left\langle K \sqcup  \mbox{\huge  $\circlearrowleft$} 
\right\rangle = 
\displaystyle\left(\frac{v^{-1}-v}{s - s^{-1}}\right)
\left\langle K \right\rangle.\]
Furthermore, for a given diagram $D=D(K)$ 
of a knot $K$,
$$\langle D \rangle  = a^{\frac{1}{2}\text{wr}(D)}\bar{\mathcal{H}}
(K;q,a)\text{\, under }s\mapsto q^{\frac{1}{2}},v\mapsto 
a^{-\frac{1}{2}},
$$
tying the variables $s$, $v$ used in \cite{HM} to the variables 
$q$, $a$ used elsewhere in this paper; $\hbox{wr}(D)$
is the {\em writhe\,} of $D$ (see there).

\subsubsection{\sf The meridian maps}
\Yboxdim4pt
Let $\varphi: \mathcal{C}\rightarrow\mathcal{C}$ be the 
\emph{meridian map} induced by adding a single oriented, unknotted 
meridian to any diagram in $S^1 \times I$ and extending linearly to 
$\mathcal{C}$.  Let $\bar{\varphi}$ be the map induced by adding a 
meridian with an orientation opposite that of $\varphi$.  Then, 
$\varphi,\bar{\varphi}$ are diagonal in their common eigenbasis 
$\{Q_{\lambda,\mu}\}\subset\mathcal{C}$ indexed by pairs 
$\lambda,\mu$ of partitions. 

The subalgebras of $\mathcal{C}$ spanned by 
$\{Q_{\lambda,\varnothing}\}$ and $\{Q_{\varnothing,\mu}\}$ are 
each isomorphic to the ring of symmetric functions in infinitely 
many variables.  Under these isomorphisms, these bases are 
identified with the basis of Schur polynomials.  Accordingly, the 
full basis $\{Q_{\lambda,\mu}\}$ is the skein-theoretic analog of 
the characters for composite partitions in \cite{K} that we 
discussed in Section \ref{Koike}. 

Now to a diagram $D$ of a knot $K$ and a composite partition 
$[\lambda,\mu]$, associate the satellite link $D\star 
Q_{\lambda,\mu}$, whose companion is $D$ and whose pattern is 
$Q_{\lambda,\mu}$.  We then have that 
$$\bar{\mathcal{H}}_{[\lambda,\mu]}(K) = 
v^{\text{wr}(D)}\langle D\star Q_{\lambda,\mu}\rangle,\
\text{wr}(D)=\text{writhe of } D,$$
i.e. the corresponding composite, unnormalized HOMFLY-PT polynomial 
for $K$ is equal to the framed, uncolored HOMFLY-PT polynomial for 
$D\star Q_{\lambda,\mu}$.

The pattern $Q_{\lambda,\mu}$ can be computed explicitly as the 
determinant of a matrix whose entries are certain idempotents 
$\{h_i,h_i^\ast\}\subset \mathcal{C}$.  For the convenience of the 
reader, some patterns for $[\lambda,\mu]$ considered in 
this paper are included in the table below. 

$$
\begin{tabular}{ |l|l| }
\hline
$[\lambda,\mu]$ & $Q_{\lambda,\mu}$\\ \hline\hline
$[\yng(1),\yng(1)]$ & $h_1h_1^\ast - 1$\\ \hline
$[\yng(1),\yng(1,1)]$ & $h_1h_1^\ast h_1^\ast - h_1h_2^\ast 
- h_1^\ast$\\ \hline
$[\yng(2),\yng(1)]$ & $h_2h_1^\ast - h_1$\\ \hline
$[\yng(1),\yng(1,1,1)]$ & $h_1h_1^\ast h_1^\ast h_1^\ast 
+ h_1h_3^\ast  + h_2^\ast - h_1h_1^\ast h_1^\ast 
- h_1h_1^\ast h_2^\ast - h_1^\ast h_1^\ast$\\ \hline
$[\yng(2,1),\yng(1)]$ & $h_1h_2h_1^\ast - h_1h_1 
- h_3h_1^\ast$\\ \hline
\end{tabular}
$$
The idempotents $h_i$ are closures of linear combinations of 
upward-oriented braids $b_i \in \Upsilon[B_i]$:
$$b_1 = 1 = \,\, \uparrow\, \in \Upsilon[B_1],\ \,
b_2 = \frac{1}{s[2]}(1+s\sigma_1) \in \Upsilon[B_2],$$
$$b_3 = \frac{1}{s^3[2][3]}(1+s\sigma_1)(1+s\sigma_2
+ s^2\sigma_2\sigma_1)\in\Upsilon[B_3],$$
in the annulus by homotopically nontrivial, 
counterclockwise-oriented strands.  Here $B_i$ is the ordinary
braid group on $i$ strands, and
the quantum integers are denoted by
$[k]\equal\frac{s^k - s^{-k}}{s - s^{-1}}$
(only in this section).  
The elements $h_i^\ast$ are then obtained by rotating the diagrams 
for $h_i$ about their horizontal axes.  That is, $h_i^\ast$ are 
linear combinations of closures of downward-oriented braids by 
clockwise-oriented strands. 

In fact, the pattern $Q_{\lambda,\mu}$ for a \emph{composite} 
partition $[\lambda,\mu]$ is distinguished by the fact that, in 
general, it contains strands oriented in both directions (clockwise 
and counterclockwise) around $S^1\times I$.  On the other hand, the 
pattern $Q_\lambda = Q_{[\lambda,\varnothing]}$ for an ordinary 
partition will consist in strands oriented all in the same 
direction. 

Let $K_{[\lambda,\mu]}\equal \frac{\langle K\star 
Q_{[\lambda,\mu]}\rangle} {\langle Q_{[\lambda,\mu]}\rangle}$, 
which is well-defined on diagrams for $K$ up to a framing 
coefficient, i.e. power of $v$.  In \cite{HM} the authors compute 
\begin{align}\label{T32HM}
K_{[\yng(1),\yng(1)]}(z,v) 
= v^2 - 4 v^4& + 4 v^6 
+ z^2 (1 + 2 v^2 - 7 v^4 + 4 v^6)\\ 
+\, &z^4 (v^2 - 2 v^4 + v^6) \for K=T^{3,2}\notag
\end{align}
in terms of variables $v$ and $z \equal s-s^{-1}$.  The 
relation to $a,q\,$ we use in this paper 
is $v=a^{-\frac{1}{2}}$ and $z = 
q^{\frac{1}{2}}-q^{-\frac{1}{2}}$;\, see below.

\subsection{\bf Rosso-Jones formula}
\subsubsection{\sf The usual theory}
The \emph{Rosso-Jones formula} \cite{JR} and its variants, e.g. 
\cite{GMV,LZ,Ste,MM}, expand the HOMFLY-PT polynomial for 
the $(\rr,\ss)$\=torus knot and a partition $\lambda\vdash n$ 
in terms of the quantum dimensions of certain irreducible 
representations:
\begin{equation}\label{RossJ}\theta_\lambda^{\rr\ss}
\bar{\mathcal{H}}_\lambda(T^{\rr,\ss}) = 
\displaystyle\sum_{\mu \hspace{1pt} 
\vdash\hspace{2pt} \rr n}c^\mu_{\lambda;\rr}
\theta^{\frac{\ss}{\rr}}_\mu\text{dim}_{q,a}(V_\mu).
\end{equation}
The formulas for $\th_\la, \theta_\mu$ and
the coefficients $c^\mu_{\lambda;\rr}$ are provided below  
in (\ref{braiding}), (\ref{Sym});\, 
$c^\mu_{\lambda;\rr}$ is nonzero only if $V_\mu$ is an irreducible 
summand of $V^{\otimes \rr}_\lambda$.  Here 
$\theta_\lambda^{\rr\ss}, \theta^{\frac{\ss}{\rr}}_\mu$
are powers, fractional for the latter.
Note that (\ref{RossJ}) 
gives the {\em unnormalized\,} polynomial as defined 
in (\ref{norm}).

\subsubsection{\sf The composite theory}
We are going to generalize the Rosso-Jones formula to the case
of composite partitions $[\la,\mu]$.  The stabilization of the 
corresponding  expansion is not\, \emph{a priori\,} clear.  
We will use the results of \cite{K} described in Section 
\ref{Koike}. The following proposition matches
formula (C.6) \cite{GJKS} (Chern-Simons theory).

\begin{prop}\label{PropRJ}
For any torus knot $T^{\rr,\ss}$ and composite partition $[\la,\mu]$ 
the corresponding (unnormalized) HOMFLY-PT polynomial admits an 
expansion:
\begin{equation}\label{CompRJ}
\theta_{[\lambda,\mu]}^{\rr\ss}\bar{\mathcal{H}}_{[\lambda,\mu]}
(T^{\rr,\ss}) = \displaystyle\sum_{[\beta,\gamma]}
c^{[\beta,\gamma]}_{[\lambda,\mu];\rr}
\theta^{\frac{\ss}{\rr}}_{[\beta,\gamma]}
\text{dim}_{q,a}(V_{[\beta,\gamma]}),
\end{equation}
into finitely many terms for which the 
$c^{[\beta,\gamma]}_{[\lambda,\mu];\rr}$ are nonzero.
Here $\theta_{[\lambda,\mu]}$ and $\theta_{[\beta,\gamma]}$ 
and the coefficients 
$c^{[\beta,\gamma]}_{[\lambda,\mu];\rr}$
are provided in (\ref{CompBE}) and (\ref{CompSym}).  
\end{prop}

{\it Proof.} First of all, it is clear from
(\ref{CompSym}) that 
$c^{[\beta,\gamma]}_{[\lambda,\mu];\rr}$ is nonzero for only 
finitely many $[\beta,\gamma]$.
Then, by construction, the 
resulting expansion (\ref{CompRJ}) will satisfy the (infinitely 
many) specializations 
\begin{equation}
\mathcal{H}_{[\lambda,\mu]}(T^{\rr,\ss};q,q^N) = 
\mathcal{H}_{[\lambda,\mu]_N}(T^{\rr,\ss};q,q^N) = 
\j^{\mathfrak{sl}_N}_{[\lambda,\mu]_N}(T^{\rr,\ss}; q),
\end{equation}
which (uniquely) define the corresponding composite HOMFLY-PT 
polynomial.  

We will divide the proof of (\ref{CompRJ})
into several intermediate steps. In what follows, 
any occurrences of $q^N$ will be replaced 
by $a$; all fractional exponents of $N$ will cancel 
in the final formula.

\subsubsection{\sf Braiding eigenvalues}
The constants $\theta_\lambda\in\mathbb{Z}[q^{\pm 1},a^{\pm 1}]$ 
in (\ref{RossJ}) are \emph{braiding eigenvalues}
from \cite{AM} , and they are 
\begin{equation}\label{braiding}
\theta_\lambda = q^{-(\kappa_\lambda + nN 
- \frac{n^2}{N})/2}\text{\, for\,  }\kappa_\lambda \equal 
\displaystyle\sum_{x\in\lambda} 2c(x),
\end{equation}
where the \emph{content} of the box $x\in \lambda$ in the 
$i${\tiny th} row and $j${\tiny th} column is $c(x) \equal j-i$.    

Now, for a composite partition $[\lambda,\mu]$ such that 
$\la\vdash m$ and $\mu\vdash n$, observe that 
$[\lambda,\mu]_N \vdash c \equal (n-m + \lambda_1N)$.  We would 
like to construct a $\kappa_{[\lambda,\mu]}$ such that 
$\kappa_{[\lambda,\mu]}|_{N=k} = \kappa_{[\lambda,\mu]_k}$ for 
any $k$.  To this end, we divide the Young diagram for 
$[\lambda,\mu]_N$ into two natural parts and count their 
individual contributions to $\kappa_{[\lambda,\mu]_N}$.
Namely, 
\begin{enumerate}
\item $\mu$ contributes $\kappa_\mu + 2\lambda_1|\mu|$
to $\kappa_{[\lambda,\mu]_N}$ and
\item $\lambda^\ast$ contributes $\kappa_{\lambda^\ast}
\!=\!\kappa_\lambda + N\lambda_1(\lambda_1+1)
-\lambda_1N(N+1)-2|\lambda|(\lambda_1-N)$.
\end{enumerate}
Thus, we can set
\begin{equation*}
\kappa_{[\lambda,\mu]}\equal\kappa_\lambda + \kappa_\mu + 
N\lambda_1(\lambda_1+1) - \lambda_1N(N+1) + 2\lambda_1|\mu| - 
2|\lambda|(\lambda_1 - N), 
\end{equation*}
so that $\kappa_{[\lambda,\mu]}|_{N=k} = 
\kappa_{[\lambda,\mu]_k}$ for any $k$, as desired.  
Furthermore we can define the composite braiding eigenvalues:
\begin{equation}\label{CompBE}
\theta_{[\lambda,\mu]}\equal q^{-(\kappa_{[\lambda,\mu]} 
+ cN - \frac{c^2}{N})/2}.
\end{equation}
One has that $\theta_{[\lambda,\mu]} 
\overset{a\mapsto q^N}{=\joinrel=\joinrel=} 
\theta_{[\lambda,\mu]_N}$ by construction.
\smallskip

The following is the key part of the proof
of Proposition \ref{PropRJ}.

\subsubsection{\sf Adams operation}
We will use Section \ref{Koike}, where we explained that the Schur 
functions $s_\lambda(\vec{x})\in\Lambda_x$ are characters for the 
irreducible polynomial representations $V_\lambda$ and described 
some of their properties.  For applications to the 
Rosso-Jones formula  we need to understand the 
$\rr$\=\emph{Adams operation} $\psi_\rr$ on $s_\lambda$;
see \cite{GMV,MM}. 

Let $p_\rr\equal\displaystyle\sum_{i\geq i} x_i^\rr\in\Lambda_x$ 
be the \emph{degree-$\rr$ power sum symmetric function}.  Then 
the $\rr$\=Adams operation on $s_\la$ may be defined formally by 
the plethysm $\psi_\rr(s_\lambda) \equal p_\rr \circ s_\lambda$.
This means that $\psi_\rr(s_\la)$ is 
determined by the coefficients $c_{\lambda; r}^\nu\in\mathbb{Z}$ 
in the expansion
\begin{equation}\label{Adams}
s_\lambda(\vec{x}^\rr) = 
\displaystyle\sum_\nu c_{\lambda; r}^\nu s_\nu(\vec{x}),
\end{equation}
where $\vec{x}^\rr\equal (x_1^\rr,x_2^\rr,x_3^\rr,\ldots)$.
The coefficients here are given an explicit description in 
\cite{LZ}:
\begin{equation}\label{Sym}
c_{\lambda; r}^\nu = \displaystyle\sum_\mu 
\frac{|C_\mu|\chi^\lambda(C_\mu)\chi^\nu(C_{r\mu})}{|\mu|},
\end{equation}
where $\chi^\lambda$ is the character of the symmetric 
group corresponding to $\lambda$, and $C_\mu$ is the conjugacy 
class corresponding to $\mu$.  

We need an analog of $\psi_\rr$
for composite partitions $[\lambda,\mu]$, which must
agree with the 
ordinary Adams operation upon specification of $N$.  Thus, we 
need to switch from (\ref{Adams}) to the expansion
\begin{equation}\label{CompAdams}
s_{[\lambda,\mu]}(\vec{x}^r,\vec{y}^r) = 
\displaystyle\sum_\nu c_{[\lambda,\mu]; r}^{[\beta,\gamma]} 
s_{[\beta,\gamma]}(\vec{x},\vec{y}),
\end{equation}
where $s_{[\lambda,\mu]}(\vec{x},\vec{y})\in 
\Lambda_x\otimes\Lambda_y$ is the universal character of \cite{K}, 
described in Section \ref{Koike}.  Applying here 
the natural projection onto $\Lambda_{N-1}$, one recovers 
the following specialization of (\ref{Adams}):
\begin{equation*}
s_{[\lambda,\mu]_N}(x_1^r,\ldots,x_{N-1}^r) = 
\displaystyle\sum_\nu c_{[\lambda,\mu]_N; r}^{[\beta,\gamma]_N} 
s_{[\beta,\gamma]_N}(x_1,\ldots,x_{N-1}).
\end{equation*}
This demonstrates that $c_{[\lambda,\mu]; r}^{[\beta,\gamma]}$
from (\ref{CompAdams}) are exactly what we need, 
i.e. this formula agrees with (\ref{Adams}) upon 
specification of $N$ and therefore
can be used for the proof of Proposition \ref{PropRJ}.

Now using (\ref{KLR}), (\ref{LRK}) and (\ref{Adams}) we obtain an 
explicit expression for these coefficients:
\begin{equation}\label{CompSym}
c_{[\lambda,\mu]; r}^{[\beta,\gamma]} = 
\displaystyle\sum_{\tau,\nu,\xi,\eta,\delta,\alpha}
(-1)^{|\tau|}N_{\nu,\tau}^\lambda N_{\tau,\xi}^\mu 
c_{\nu ; r}^\eta c_{\xi ; r}^\delta N_{\beta,\alpha}^\eta 
N_{\gamma,\alpha}^\delta,
\end{equation}
where the sum is over arbitrary sextuples of Young diagrams.
Recall that $N_{\nu,\tau}^\lambda$, are the 
Littlewood-Richardson coefficients from (\ref{LR}).

Although this formula appears rather complicated, observe 
that the terms are only nonzero for relatively few (and finitely 
many) choices of $(\tau,\nu,\xi,\eta,\delta,\alpha)$.  
In light of (\ref{Sym}) and the combinatorial nature 
of the Littlewood-Richardson rule, these formula provides a 
completely combinatorial description of 
$c_{[\lambda,\mu]; r}^{[\beta,\gamma]}$.

The following is the last step of the proof.

\subsubsection{\sf Quantum dimensions}\label{qdim}
We define the $q,a$\=\emph{integer\,} by 
\begin{equation*}
[uN+v]_{q,a} \equal \frac{a^{\frac{u}{2}}q^{\frac{v}{2}} 
- a^{-\frac{u}{2}}q^{-\frac{v}{2}}}{q^{\frac{1}{2}} 
- q^{-\frac{1}{2}}},
\end{equation*}
for $u,v\in\mathbb{Z}$, where $N$ is ``generic", i.e. it
is treated here as a formal variable. Setting here $a=q^N$ 
for $N\in \N$, we obtain the ordinary quantum integer 
$[uN+v]_q$. We will suppress the subscript ``$q,a$" in this and
the next subsection, simply writing 
$[\hspace{2pt}\cdot\hspace{2pt}]$. 

For an irreducible representation $V_\mu$, its 
\emph{stable quantum dimension\,} 
is given by the quantum Weyl dimension formula 
\begin{equation}\label{murhoal}
\text{dim}_{q,a}(V_\mu) = \displaystyle\prod_{\alpha\in A_{N-1}^+} 
\frac{[(\mu + \rho,\alpha)]}{[(\rho,\alpha)]},
\end{equation}
where the Young diagram $\mu$ is interpreted in the usual way as 
a weight for  $\mathfrak{sl}_N$ for generic $N$ and  
$\rho = \frac{1}{2}\sum_{\alpha>0}\alpha$
for $A_{N-1}$.  

Then it only depends on
the diagram $\mu$, which includes the actual number of
factors due to the cancelations. We note that
such a stabilization holds in the theory of Macdonald polynomials
of type $A_{N-1}$ as well; see formula (\ref{macdsym}) and
Theorem \ref{STABILIZ}, ($i$).
\smallskip
  
The stable quantum dimension for  a composite partition 
$[\beta,\gamma]$ is defined as follows:
\begin{equation}\label{CompWeyl}
\text{dim}_{q,a}(V_{[\beta,\gamma]}) \equal 
\displaystyle\prod_{\alpha\in A_{N-1}^+} \frac{[([\beta,\gamma]_N 
+ \rho,\alpha)]}{[(\rho,\alpha)]}.
\end{equation}
Similarly to (\ref{murhoal}), we claim that there is no actual 
dependence of $N$ in this formula (including the actual number
of factors). However the justification is somewhat
more involved because the weight

$$[\beta,\gamma]_N = 
\displaystyle\sum_{j=1}^{\ell(\gamma)}(\gamma_i-\gamma_{i+1})
\omega_i + \sum_{j=1}^{\ell(\beta)}(\beta_j - \beta_{j+1})
\omega_{N-j},$$
depends on $N$ (in contrast to the case of one diagram).
We will omit a straightforward justification; see 
table (\ref{Ex2}) below and the general formula (C.3) from 
\cite{GJKS} (a calculation of normalized open-string stretched
annulus amplitudes). 
Finally, the relation
$\text{dim}_{q,a}(V_{[\beta,\gamma]})|
_{a\mapsto q^N}$ $= \text{dim}_{q}(V_{[\beta,\gamma]_N})$
concludes the proof of Proposition \ref{PropRJ}. \sq  
\smallskip

Formula (\ref{CompRJ}) provides a purely combinatorial 
and computationally effective way of producing HOMFLY-PT 
polynomials for arbitrary torus knots and composite 
representations. See examples below and also Section C
from \cite{GJKS}.

\subsubsection{\sf Simplest examples}\label{Exs}
First, we evaluate the 
(ordinary) Rosso-Jones formula (\ref{RossJ}) for the trefoil 
$T^{3,2}$ and $\lambda = \square$.  The necessary values are 
contained in table (\ref{Ex1}): 
\begin{equation}\label{Ex1}
\begin{tabular}{ |l|c|c|c| }
\hline
$\mu$ & $\theta_{\mu}
$ & $c_{\yng(1);2}^{\mu}$ 
& $\text{dim}_{q,a}(V_{\mu})$\\ \hline\hline
$\yng(1)$ & $a^{-\frac{1}{2}}q^{\frac{1}{2N}}$ & $0$ & $[N]$\\ 
\hline
$\yng(2)$ & $a^{-1}q^{\frac{2}{N}-1}$ & $1$ & 
$ \frac{[N][N+1]}{[2]}$\\ \hline
$\yng(1,1)$ & $a^{-1}q^{\frac{2}{N}+1}$ & $-1$ & 
$\frac{[N-1][N]}{[2]}$\\ \hline
\end{tabular}
\end{equation}
\vskip -0.25cm 
\centerline{ \hspace{170pt}.}
\noindent
Inserting the components of (\ref{Ex1}) into formula (\ref{RossJ}), 
we obtain the familiar expression:
\begin{align*}
\mathcal{H}_{\yng(1)}(T^{3,2};q,a) \;& = \;  
\frac{ \theta_{\yng(1)}^{-6}(\theta_{\yng(2)}^{\frac{3}{2}}
\text{dim}_{q,a}(V_{\yng(2)}) - \theta_{\yng(1,1)}^{\frac{3}{2}}
\text{dim}_{q,a}(V_{\yng(1,1)}))}{\text{dim}_{q,a}(V_{\yng(1)})}\\
& = \; aq^{-1} - a^{2} + aq,
\end{align*}
the normalized HOMFLY-PT polynomial of $T^{3,2}$.  Note that 
although $\square$ appears with coefficient $0$ in the expansion 
(\ref{RossJ}), we include it in table (\ref{Ex1}) since both 
$\theta_{\yng(1)}$ and $\text{dim}_{q,a}(V_{\yng(1)})$ are needed 
to give the final, normalized polynomial, as defined in 
(\ref{norm}).

Similarly, we evaluate our composite Rosso-Jones formula 
(\ref{CompRJ}) for the trefoil $T^{3,2}$ and $[\yng(1),\yng(1)]$ 
using table (\ref{Ex2}): 
\begin{equation}\label{Ex2}
\begin{tabular}{ |l|c|c|c| }
\hline
$[\beta,\gamma]$ & $\theta_{[\beta,\gamma]}
$ & $c_{[\yng(1),\yng(1)];2}^{[\beta,\gamma]}$ 
& $\text{dim}_{q,a}(V_{[\beta,\gamma]})$\\ \hline\hline
$[\yng(1),\yng(1)]$ & $a^{-1}$ & $0$ & $[N-1][N+1]$\\ \hline
$[\yng(2),\yng(2)]$ & $q^{-2}a^{-2}$ & $1$ & 
$\frac{[N-1][N]^2[N+3]}{[2][2]}$\\ \hline
$[\yng(2),\yng(1,1)]$ & $a^{-2}$ & $-1$ & 
$\frac{[N-2][N-1][N+1][N+2]}{[2][2]}$\\ \hline
$[\yng(1,1),\yng(2)]$ & $a^{-2}$ & $-1$ & 
$\frac{[N-2][N-1][N+1][N+2]}{[2][2]}$\\ \hline
$[\yng(1,1),\yng(1,1)]$ & $q^2a^{-2}$ & 
$1$ & $\frac{[N-3][N]^2[N+1]}{[2][2]}$\\ \hline
$[\varnothing,\varnothing]$ & $1$ & $1$ & $1$\\ \hline
\end{tabular}
\end{equation}
\vskip -0.25cm 
\centerline{ \hspace{235pt}.}
\noindent
Inserting the components of (\ref{Ex2}) into formula 
(\ref{CompRJ}), we obtain
\begin{align*}
\mathcal{H}_{[\yng(1),\yng(1)]}(T^{3,2};q,a) \; = \; 
& a^2(q^{-2} + q^2 + 2)
+ a^3(-2q^{-2} + q^{-1} + q - 2q^2 - 2)\\
& + a^4(q^{-2} - 2q^{-1} - 2q + q^2 + 3)
+ a^5(q^{-1} + q - 2),
\end{align*}
where we include $[\yng(1),\yng(1)]$ in table (\ref{Ex2}) for the 
same reason that we included $\square$ in table (\ref{Ex1}).

Observe that we can touch base with formula 
(\ref{T32HM}) from \cite{HM} by 
$$a^5T^{3,2}_{[\yng(1),\yng(1)]}(q^{\frac{1}{2}} 
- q^{-\frac{1}{2}}, a^{-\frac{1}{2}}) 
= \mathcal{H}_{[\yng(1),\yng(1)]}(T^{3,2};q,a).$$
Our expression for 
$\mathcal{H}_{[\yng(1),\yng(1)]}(T^{3,2};q,a)$ agrees with that 
obtained in \cite{PBR}. See also examples (C.8-16) from
\cite{GJKS}, obtained there via Chern-Simons theory 
(open-string amplitudes); they match our ones.

\setcounter{equation}{0}
\section{\sc DAHA superpolynomials}\label{DAHA}
\subsection{\bf Definition of DAHA}
\subsubsection{\sf Affine root systems}
Let $R=\{\al\}   \subset \R^n$ be a root system of type
$A_n,\ldots,\!G_2$ with respect to a euclidean form 
$(\,,\,)$ on $\R^n$, normalized by the condition  
$(\al,\al)=2$ for {\em short\,} roots.
Let $W=\lan s_\al\ran$ be its Weyl group, and let
$R_{+}$ be the set of positive roots
corresponding to a fixed set $\{\al_1,...,\al_n\}$ 
of simple roots for $R$.
The weight lattice is
$P=\oplus^n_{i=1}\Z \om_i$, 
where $\{\om_i\}$ are fundamental weights:
$ (\om_i,\al_j^\vee)=\de_{ij}$ for the
coroots $\al^\vee=2\al/(\al,\al)$;
$P_{\pm}=\oplus^n_{i=1}\Z_{\pm} \om_i$, 
for $\Z_{\pm}=\{m\in\Z, \pm m\ge 0\}$. 

Setting 
$\nu_\al\equal (\al,\al)/2$,
the vectors $\ \tal=[\al,\nu_\al j] \in
\R^n\times \R \subset \R^{n+1}$
for $\al \in R, j \in \Z $ form the
{\em twisted affine root system\,}
$\tR \supset R$ ($z\in \R^n$ are identified with $ [z,0]$).
We add $\al_0 \equal [-\vth,1]$ to the simple
roots for the {\em maximal short root\,} $\vth\in R_+$.
The corresponding set
$\tR_+$ of positive roots is 
$R_+\cup \{[\al,\nu_\al j],\ \al\in R, \ j > 0\}$.

The set of the indices of the images of $\al_0$ by all
automorphisms of the affine Dynkin diagram will be denoted by 
$O$; let $O'\equal\{r\in O, r\neq 0\}$.
The elements $\om_r$ for $r\in O'$ are  
{\em minuscule weights\,}. We set $\om_0=0$.
\smallskip

\subsubsection{}{\sf Extended Weyl group.}
Given $\tal=[\al,\nu_\al j]\in \tR,  \ b \in P$, let
\begin{align}
&s_{\tal}(\tz)\ =\  \tz-(z,\al^\vee)\tal,\
\ b'(\tz)\ =\ [z,\ze-(z,b)]
\label{ondon}
\end{align}
for $\tz=[z,\ze] \in \R^{n+1}$.
The
{\em affine Weyl group\,} $\tW=\lan s_{\tal}, \tal\in \tR_+\ran$ 
is the semidirect product $W\lsmash Q$ of
its subgroups $W=$ $\lan s_\al,
\al \in R_+\ran$ and $Q$, where $\al$ is identified with
\begin{align*}
& s_{\al}s_{[\al,\,\nu_{\al}]}=\
s_{[-\al,\,\nu_\al]}s_{\al}\for
\al\in R.
\end{align*}

The {\em extended Weyl group\,} $ \hW$ is $W\lsmash P$, where
the corresponding action is 
\begin{align}
&(wb)([z,\ze])\ =\ [w(z),\ze-(z,b)] \for w\in W, b\in P.
\label{ondthr}
\end{align}
It is isomorphic to $\tW\lsmash \Pi$ for $\Pi\equal P/Q$. 
The latter group consists of $\pi_0=$id\, and the images $\pi_r$
of minuscule $\om_r$ in $P/Q$. 

The group $\Pi$
is naturally identified with the subgroup of $\hW$ of the
elements of the length zero; the {\em length\,} is defined as 
follows:
\begin{align*}
&l(\hw)=|\la(\hw)| \for \la(\hw)\equal\tR_+\cap \hw^{-1}(-\tR_+).
\end{align*}
One has $\om_r=\pi_r u_r$ for $r\in O'$, where $u_r$ is the 
element $u\in W$ of minimal length such that $u(\om_r)\in P_-$. 
\smallskip

Setting $\hw = \pi_r\tw \in \hW$ for $\pi_r\in \Pi,\, \tw\in \tW,$
\,$l(\hw)$ coincides with the length of any reduced decomposition
of $\tw$ in terms of the simple reflections
$s_i,$ $0\le i\le n.$ 

\subsubsection{\sf Parameters} 
We follow \cite{CJJ,CJ,C101}.
Let $\mm,$ be the least natural number
such that  $(P,P)=(1/\mm)\Z.$  Thus
$\mm=|\Pi|$ unless 
$\mm=2 \for D_{2k}$ and $\ \mm=1 \for B_{2k},C_{k}.$ 

The double affine Hecke algebra, {\em DAHA\,}, depends
on the parameters
$q, t_\nu\, (\nu\in \{\nu_\al\})\,$ and is naturally defined
over the ring
$Z_{q,t}\equal\Z[q^{\pm 1/\mm},t_\nu^{\pm 1/2}]$
formed by
polynomials in terms of $q^{\pm 1/\mm}$ and
$\{t_\nu^{1/2}\}.$ 

For $\tal=[\al,\nu_\al j] \in \tR,\ 0\le i\le n$, we set
\begin{align*}
&   t_{\tal} =t_{\al}=t_{\nu_\al}=q_\al^{k_\nu} ,\ \, 
q_{\tal}=q^{\nu_\al}, \ \, t_i = t_{\al_i},\ \,q_i=q_{\al_i},
\end{align*}

Also, using here (and below) {\em\small\, sht,\ lng\,} instead 
of $\nu$, we set
\begin{align*}
\rho_k\equal \frac{1}{2}\!\sum_{\al>0} k_\al \al=
k_{\sht}\rho_{\sht}\!+\!k_{\lng}\rho_{\lng},\ \,
\rho_\nu=\frac{1}{2}\!\sum_{\nu_\al=\nu} \al=
\!\!\sum_{\nu_i=\nu,\, i>0}  \om_i.
\end{align*}

For pairwise commutative $X_1,\ldots,X_n,$
\begin{align}
& X_{\tb}\ \equal\ \prod_{i=1}^nX_i^{l_i} q^{ j}
\iif \tb=[b,j],\ \hw(X_{\tb})\ =\ X_{\hw(\tb)},
\label{Xdex}\\
&\hbox{where\ } b=\sum_{i=1}^n l_i \om_i\in P,\ j \in
\frac{1}{ m}\Z,\ \hw\in \hW.
\notag \end{align}
For instance, $X_0\equal X_{\al_0}=qX_\vth^{-1}$.
\medskip

\subsubsection{\sf The main definition}
Recall that 
$\om_r=\pi_r u_r$ for $r\in O'$ (see above).
We will use that $\pi_r^{-1}$ is $\pi_{\iota(i)}$, where
$\iota$ is the standard involution  of the nonaffine 
Dynkin diagram,
induced by $\al_i\mapsto -w_0(\al_i)$. Generally,
$\iota(b)=-w_0(b)=b^\iota$, where $w_0$ is the longest element 
in $W$.  Finally, we set $m_{ij}=2,3,4,6$
when the number of links between $\al_i$ and $\al_j$ in the affine 
Dynkin diagram is $0,1,2,3$.

\begin{definition}
The double affine Hecke algebra $\HH\ $
is generated over $\Z_{q,t}$ by
the elements $\{ T_i,\ 0\le i\le n\}$,
pairwise commutative $\{X_b, \ b\in P\}$ satisfying
(\ref{Xdex})
and the group $\Pi,$ where the following relations are imposed:

(o)\ \  $ (T_i-t_i^{1/2})(T_i+t_i^{-1/2})\ =\
0,\ 0\ \le\ i\ \le\ n$;

(i)\ \ \ $ T_iT_jT_i...\ =\ T_jT_iT_j...,\ m_{ij}$
factors on each side;

(ii)\ \   $ \pi_rT_i\pi_r^{-1}\ =\ T_j \iif
\pi_r(\al_i)=\al_j$;

(iii)\  $T_iX_b \ =\ X_b X_{\al_i}^{-1} T_i^{-1} \iif
(b,\al^\vee_i)=1,\
0 \le i\le  n$;

(iv)\ $T_iX_b\ =\ X_b T_i\ $ if $\ (b,\al^\vee_i)=0
\for 0 \le i\le  n$;

(v)\ \ $\pi_rX_b \pi_r^{-1}\ =\ X_{\pi_r(b)}\ =\
X_{ u^{-1}_r(b)}
 q^{(\om_{\iota(r)},b)},\  r\in O'$.
\label{double}
\end{definition}

Given $\tw \in \tW, r\in O,\ $ the product
\begin{align}
&T_{\pi_r\tw}\equal \pi_r T_{i_l}\cdots T_{i_1},\where
\tw=s_{i_l}\cdots s_{i_1} \for l=l(\tw),
\label{Twx}
\end{align}
does not depend on the choice of the reduced decomposition.
Moreover,
\begin{align}
&T_{\hv}T_{\hw}\ =\ T_{\hv\hw}\  \hbox{ whenever\,}\
 l(\hv\hw)=l(\hv)+l(\hw) \for
\hv,\hw \in \hW. \label{TTx}
\end{align}
In particular, we arrive at the pairwise
commutative elements 
\begin{align}
& Y_{b}\equal
\prod_{i=1}^nY_i^{l_i} \iif
b=\sum_{i=1}^n l_i\om_i\in P,\ 
Y_i\equal T_{\om_i},b\in P.
\label{Ybx}
\end{align}
When acting in the polynomial representation, they are 
called {\em difference Dunkl operators.}
\smallskip

\subsubsection{\sf Automorphisms}
The following maps can be (uniquely) extended to
an automorphism of $\HH\,$, fixing $\ t_\nu,\ q$
and their fractional powers;
see \cite{C101}, (3.2.10)-(3.2.15). Adding 
$q^{1/(2\mm)}$ to $\Z_{q,t}$,
\begin{align}\label{tauplus}
& \tau_+:\  X_b \mapsto X_b, \ T_i\mapsto T_i\, (i>0),\
\ Y_r \mapsto X_rY_r q^{-\frac{(\om_r,\om_r)}{2}}\,,
\\
& \tau_+:\ T_0\mapsto  q^{-1}\,X_\vth T_0^{-1},\
\pi_r \mapsto q^{-\frac{(\om_r,\om_r)}{2}}X_r\pi_r\
(r\in O'),\notag\\
& \label{taumin}
\tau_-:\ Y_b \mapsto \,Y_b, \ T_i\mapsto T_i\, (i\ge 0),\
\ X_r \mapsto Y_r X_r q^\frac{(\om_r,\om_r)}{ 2},\\
&\tau_-(X_{\vth})= 
q T_0 X_\vth^{-1} T_{s_{\vth}}^{-1};\ \
\si\equal \tau_+\tau_-^{-1}\tau_+\, =\,
\tau_-^{-1}\tau_+\tau_-^{-1},\notag\\
&\si(X_b)=Y_b^{-1},\   \si(Y_b)=
T_{w_0}^{-1}X_{b^\iota}^{-1}T_{w_0},\ \si(T_i)=T_i (i>0).
\label{taux}
\end{align}

The group  $PSL_{\,2}^\wedge(\Z)\,$  generated
by $\tau_{\pm}$, the {\em projective
$PSL_2(\Z)$\,} due to Steinberg, has a  natural projection 
onto $PSL_2(\Z)$, corresponding to taking
$t^{1/(2\mm)}_\nu=1=q^{1/(2\mm)}$:\ \  
$\tau_+\mapsto$ 
{\tiny 
$\begin{pmatrix}1 & 1 \\0 & 1 \\ \end{pmatrix}$},\ 
$\tau_-\mapsto$ 
{\tiny 
$\begin{pmatrix}1 & 0 \\1 & 1 \\ \end{pmatrix}$},\
$\si\mapsto$ 
{\tiny 
$\begin{pmatrix}0 & 1 \\-1 & 0 \\ \end{pmatrix}$}.\

\subsection{\bf DAHA-Jones polynomials}
\subsubsection{\sf Coinvariant}
Following \cite{C101}, 
we use the PBW Theorem to express any $H\in \HH$ in the form 
\,$\sum_{b,w,c} d_{b,w,c}\, X_b T_{w} Y_c$\, for $w\in W$,
$b,c\in P$ (this presentation is unique). Then we substitute:
\begin{align}\label{evfunct}
\{\,\}_{ev}:\ X_b \ \mapsto\  q^{-(\rho_k,b)},\ 
Y_c \ \mapsto\  q^{(\rho_k,c)},\ 
T_i \ \mapsto\  t_i^{1/2}. 
\end{align}

The functional $\,\HH\ni H\mapsto \{H\}_{ev}$, 
called {\em coinvariant\,}, acts via the projection 
$H\mapsto H(1)$ of $\HH\,$
onto the {\em polynomial representation \,}$\v$, which is
the $\HH$\=module induced from the one-dimensional
character $T_i(1)=t_i^{-1/2}=Y_i(1)$ for $1\le i\le n$ and
$T_0(1)=t_0^{-1/2}$. Recall  that $t_0=t_{\sht}$; 
see \cite{C101,CJ}.
\smallskip

\subsubsection{\sf Macdonald polynomials}
The polynomial representation
is isomorphic to $\Z_{q,t}[X_b]$
as a vector space, and the action of $T_i(0\le i\le n)$ there 
is given by 
the {\em Demazure-Lusztig operators\,}:
\begin{align}
&T_i\  = \  t_i^{1/2} s_i\ +\
(t_i^{1/2}-t_i^{-1/2})(X_{\al_i}-1)^{-1}(s_i-1),
\ 0\le i\le n.
\label{Demazx}
\end{align}
The elements $X_b$ become the multiplication operators 
and  $\pi_r (r\in O')$ act via the general formula
$\hw(X_b)=X_{\hw(b)}$ for $\hw\in \hW$. 

The Macdonald polynomials $P_b(X)$ are uniquely defined
as follows. Let $c_+$ be the unique element 
such that $c_+\in W(c)\cap P_+$. For $b\in P_+$,
\begin{align*}
&P_b\! -\!\!\!\!\sum_{b'\in W(b)}\!\!\! X_{b'} 
\in\, \oplus_{b_+\neq c_+\in b+Q_+}\Q(q,t_\nu) X_c 
\hbox{\, and\, }
CT\bigl(P_b X_{c^\iota}\,\mu(X;q,t)\bigr)\!=\!0
\\
&\hbox{for such $c$,\, where\, } 
\mu(X;q,t)\equal\!\prod_{\al \in R_+}
\prod_{j=0}^\infty \frac{(1\!-\!X_\al q_\al^{j})
(1\!-\!X_\al^{-1}q_\al^{j+1})
}{
(1\!-\!X_\al t_\al q_\al^{j})
(1\!-\!X_\al^{-1}t_\al^{}q_\al^{j+1})}\,.
\end{align*}
Here $CT$ is the constant term; 
$\mu$ is considered
a Laurent series in $X_b$ with 
the coefficients expanded in terms of
positive powers of $q$. The coefficients of
$P_b$ belong to the field $\Q(q,t_\nu)$.
One has:
\begin{align}\label{macdsym}
&P_b(X^{-1})\,=\,P_{b^\iota}(X)\,=\,
P_{b}(q^{-\rho_k})=P_{b}(q^{\rho_k})\\
\label{macdeval}
&\, =\,q^{-(\rho_k,b)}\,
\prod_{\al>0}\,\prod_{j=0}^{(\al^{\!\vee},b)-1}
\,\Bigl(
\frac{
1- q_\al^{j}t_\al X_\al(q^{\rho_k})
 }{
1- q_\al^{j}X_\al(q^{\rho_k})
}
\Bigr).
\end{align}
See \cite{C101}, formula (3.3.23); recall that 
$\iota(b)=b^\iota=-w_0(b)$ for $b\in P$.

\subsubsection{\sf DAHA-Jones polynomials}
We begin with the following theorem, which is
from \cite{CJ,CJJ}.


Torus knots $T^{\rr,\ss}$ are naturally represented   
by $\ga_{\rr,\ss}\in PSL_{\,2}(\Z)$ with the 
first column $(\rr,\ss)^{tr}$ ($tr$ is the transposition)
for $\,\rr,\ss\in \N$, assuming that \,gcd$(\rr,\ss)=1$. 
Let $\hat{\ga}_{\rr,\ss}\in PSL_{\,2}^{\wedge}(\Z)$ be
any pullback of $\ga_{\rr,\ss}$.
\smallskip

For a polynomial $F$ in terms of 
fractional powers of $q$ and $t_\nu$, 
the {\em tilde-normalization}
$\tilde{F}$ will be the result of the division of $F$
by the lowest $q,t_\nu$\=monomial, assuming that it 
is well defined. We put $q^\bullet t^\bullet$ for
a monomial factor (possibly fractional)
in terms of $q,t_\nu$.

\begin{theorem}\label{JONITER}
Given a torus knot $T^{\rr,\ss}$, 
we lift $(\rr,\ss)^{tr}$ to
$\ga$ and then to $\hat{\ga}\in PSL_{\,2}^{\wedge}(\Z)$ as above.

(i) The {\sf DAHA-Jones (or refined) polynomial} for a reduced
irreducible root 
system $R$ and $b\in P_+\,$ is defined as follows:
\begin{align}\label{jones-dit}
&J\!D_{\,\rr,\,\ss\,}^{R}(b\,;\,q,t)\ =\ 
J\!D_{\,\rr,\,\ss\,}(b\,;\,q,t)\, \equal\,
\bigl\{\hat{\ga}(P_b)\bigr\}_{ev}.
\end{align}

(ii) It does not depend on the ordering of $\rr,\ss$ or
on the particular choice of \,$\ga\in PSL_2(\Z),$
$\hat{\ga}\in PSL_{\,2}^{\wedge}(\Z)$. 
The tilde-normalization 
$\tilde{J\!D}_{\rr,\ss}\,(b\,;\,q,t)$ is well defined
and is a polynomial in terms of $\,q,t_\nu$ with
constant term $1$. 

(iii) {\sf Specialization at the trivial center charge.\,} For 
$b=\hbox{\small$\sum$}_{i=1}^n b_i \om_i$,
\begin{align}\label{jones-eval}
J\!D_{\,\rr,\,\ss}
\,\bigl(b\,;\,q\!=\!1,t
\bigr)\!=\!
\hbox{\small$\prod$}_{i=1}^n J\!D_{\,\rr,\,\ss}\,
(\om_i\,;\,q\!=\!1,t)^{b_i} \hbox{\, for any \,} \rr,\ss.
\end{align}
\vskip -1.25cm \sq
 \end{theorem}

It was conjectured in \cite{CJ} in general (and checked there 
for $A_n$) that  $J\!D_{\,\rr,\,\ss}\,(b\,;\,q, t_\nu
\!\mapsto\! q_\nu)$ 
coincide up to $q^\bullet$ 
with the reduced Quantum Group (WRT) invariants for the
corresponding  $T^{\rr,\ss}$ and any colors $b\in P_+$.
The Quantum Group is associated with the {\em twisted\,}
root system $\tR$. The {\em shift operator\,} was used there
to deduce this coincidence from \cite{LZ,Ste} in the case
of $A_n$ and torus knots. The papers \cite{Ste,CC} provide
the necessary tools to establish this coincidence for $D_n$.
Quite a few further confirmations for other root systems 
are known by now; the second author (R.~E.) checked such a  
coincidence with the DAHA formulas provided
in \cite{CJ} for the minuscule and quasi-minuscule weights 
for $E_6$ (unpublished).


\subsection{\bf DAHA superpolynomials}
Theorem \ref{JONITER} leads to the theory of 
{\em DAHA-superpolynomials\,}, 
which are the result of the {\em stabilization\,}  of 
$\tilde{J\!D}^{A_n}(b;q,t)$ with respect to $n$.
This stabilization was announced 
in \cite{CJ}; its proof was published in \cite{GoN}. 
Both approaches use \cite{SV}; we note that the stabilization 
holds for arbitrary torus iterated knots. 

Following \cite{SV} (see also \cite{GoN,CJJ}),
we can generalize the stabilization
construction to the torus knots in the annulus.
\smallskip

The pairs $\{\rr,\ss\}$ remains the same,
but now colored torus knots $T^{\rr,\ss}$ will 
be treated as link diagrams in the annulus; see
Section \ref{HOMFLY}.

\begin{theorem}\label{STABILIZ}
We switch to $A_n$, setting $t=t_{\sht}=q^k$. 
Let $b,c\in P_+^n\equal P_+^{A_n}$ and  
$\la,\mu$ be the corresponding
Young diagrams (with no greater than $n$ rows). Recall that
$[\la,\mu]_N\in P_+^{N-1}$
is $\,b^*+c$, where $\,N\ge \ell(\la)+\ell(\mu)\,$ 
and $\,(\om_i)^\ast= \om_{N-i}$\,;
see (\ref{lamuN}).

(i) {\sf Stabilization.} Given a pair $\{\rr,\ss\}$
as above, 
there exists a polynomial 
$H\!D_{\,\rr,\,\ss}\,([\la,\mu]\,;\,q,t,a)$
from $\Z[q,t^{\pm 1},a]$ 
such that its coefficient of $a^0$ is tilde-normalized 
(i.e. in the form $\sum_{u,v\ge 0}C_{u,v}q^u t^v$ 
with $C_{0,0}=1$) and 
\begin{align}\label{jones-sup}
H\!D_{\rr,\ss}([\la,\mu]\,;q,t,a\!=\!-t^{N}) =
\tilde{J\!D}_{\rr,\ss}^{A_{N-1}}(b^*\!+c\,;q,t)
\hbox{\, for any\, } N\!>\! n.
\end{align} 
This polynomial does not depend on the ordering of $\rr,\ss$
or that of $\la,\mu$.

(ii) {\sf Specialization at $q=1$.}
Setting 
$H\!D_{\,\rr,\,\ss}\,(\la)=H\!D_{\,\rr,\,\ss}\,([\varnothing,\la])$,
\begin{align}\label{jones-seval}
H\!D_{\,\rr,\,\ss}\,([\la,\mu]\,;&\,q\!=\!1,t,a)\\
=H\!D_{\,\rr,\,\ss}\,&(\la\,;q\!=\!1,t,a)\,
H\!D_{\,\rr,\,\ss}\,(\mu\,;q\!=\!1,t,a),\where\notag\\
H\!D_{\,\rr,\,\ss}\,(\la\,;\,q\!=\!1,t,a)\,=\,&\prod_{i=1}^n
H\!D_{\,\rr,\,\ss}\,(\om_i\,;\,q\!=\!1,t,a)^{b_i} 
\hbox{\,\, for\,\, }
b\!=\!\sum_{i=1}^n b_i\om_i,\notag
\end{align}
$\,b\,$ corresponds to $\,\la\,$
and $\,\om_i\,$ means the column with $i$ boxes.
\sq
\end{theorem}

\subsubsection{\sf Degree of \!{\mathversion{bold}$a$}
and duality}
Assuming  that $\rr>\ss$, we conjecture that
\begin{align}\label{degaconj}
\hbox{deg}_a H\!D_{\rr,\,\ss}\,([\la,\mu]\,;\,q,t,a)\,
=\,\ss(|\la|+|\mu|)-|\la\!\vee\!\mu|,
\end{align}
where $\la\!\vee\!\mu$ (the join operation) is the smallest 
Young diagram containing them,
$|\la|$ is the number of boxes in $\la$.
This is based on the numerical evidence and on a 
generalization of the construction from \cite{GoN} to
the composite case (though we did not check all details).

Let us generalize the DAHA-duality from \cite{CJ}
(justified in \cite{GoN}) to the composite case; 
see also \cite{GS, CJJ}.

\begin{theorem}\label{SUPRDUAL}
{\sf Composite super-duality.} 
Up to a power of $q$ and $t$,
\begin{align}\label{comp-duality}
H\!D_{\,\rr,\,\ss}\,([\la,\mu]\,;q,t,a)=q^{\bullet}t^{\bullet} 
H\!D_{\,\rr,\,\ss}\,([\la^{tr},\mu^{tr}]\,;t^{-1},q^{-1},a), 
\end{align}
where $\la^{tr}$ is the transposition of $\la$. 
\end{theorem} 
{\it Proof.} According to the remark 
after the super-duality formula (1.44) from 
Section 1.6 of \cite{CJJ}, the standard type $A$ (one-diagram) 
duality is equivalent to $q^\bullet$-proportionality between 
$\j_{\rr,\ss}^{A_n}(\la\,;q,t)$ and  
$\j_{\rr,\ss}^{A_m}(\la^{tr}\,;t^{-1},q^{-1})$
for $t=q^{-(m+1)/(n+1)}$ (i.e. for $k=-\frac{m+1}{n+1}$)
and all possible relatively prime $m+1,n+1\in \N$. This
is directly connected with the generalized {\em level-rank
duality\,}. Using that $q, n, m$ are essentially
arbitrary,  we conclude that these proportionality
conditions (all of them) are equivalent to the duality. 
The latter was proved in \cite{GoN}; the above argument 
(and the theory of perfect DAHA modules at roots of unity
from \cite{C101}) can be used for the justification of
the standard super-duality as well (unpublished).

This reformulation of the super-duality in
terms of the DAHA-Jones polynomials (i.e. without $a$)
gives the composite super-duality upon
considering the diagrams in the form $[\la,\mu]_N$. 
\sq

Combining the evaluation formula (\ref{jones-seval})
with the duality:
\begin{align}\label{jones-sevalt}
H\!D_{\rr,\ss}\,([\la,\mu]\,;\,q,\,&t\!=\!1,a)\\
&=H\!D_{\rr,\ss}\,(\la\,;\,q,t\!=\!1,a)\,
H\!D_{\rr,\ss}\,(\mu\,;\,q,t\!=\!1,a).\notag
\end{align}

\subsubsection{\sf Color exchange}
The following theorem can be proved following
Sections 1.6, 1.7 from \cite{CJJ}.

\begin{theorem}\label{COLOREXCH}
{\sf Color Exchange.} Let $t=q^k$ for $k\in -\Q_+$.
For $\la,\mu$ as above, we assume the existence of
permutations $v,w\in \S_n$ satisfy the following conditions.
Setting $\la=\{l_1\ge l_2\ldots\ge l_n\ge 0\}$,
\begin{align}\label{lamu}
\la'=\{l'_1,\ldots,l_n'\}\equal
\{l_{u(i)}+k(i-u(i)),\ i=1,2,\ldots,n\}
\end{align} 
must be a diagram, i.e we require
$l'_i\geq l_{i+1}'$ and $l'_i\in \Z_+$. 
Similarly, $\mu'$ defined by $\mu,w$ (for the same $k$)
is assumed a Young diagram. Then 
$H\!D_{\,\rr,\,\ss}\,([\la,\mu]\,;\,q,t,a) = 
H\!D_{\,\rr,\,\ss}\,([\la',\mu']\,;\,q,t,a)$
for such $q,t$ and any $\rr,\ss$.
\end{theorem}
\vskip -0.9cm\sq

Let us provide an example for $t=q^{-\kappa}$, $\kappa\in \N$
(see \cite{CJJ}, formula (1.47) for details). For any $p>0$
and $i\in \{1,2\}\ni j$, one has:
\begin{align*}
&H\!D_{\rr,\ss}(\,[\kappa b^{(i)},\kappa b^{(j)}]
\,;\,q,q^{-\kappa},a)= q^{\bullet}\,
H\!D_{\rr,\ss}(\,[\kappa c^{(i)},\kappa c^{(j)}]\,;
\,q,q^{-\kappa},a)
\hbox{\,\, for}\\
&b^{(1)}=\om_{p+1}, c^{(1)}=(p+1)\om_1
\hbox{\, and \,} b^{(2)}=p\om_{p+1}, c^{(2)}=(p+1)\om_p\,,
\end{align*}
where the weights are identified with the corresponding
diagrams. If $\kappa=1$, then $t=q^{-1}$ and these relations
follow from the duality.

\subsubsection{\sf Obtaining HOMFLY-PT polynomials}
\begin{theorem}\label{CONNEC}{\sf HOMFLY-PT via DAHA.}
For $\rr,\ss$ and $\la,\mu$ as above,
\begin{align}\label{conthm}
&H\!D_{\rr,\ss}([\lambda,\mu]; q,t\mapsto q,a\mapsto -a) 
= \mathcal{H}_{[\lambda,\mu]}(T^{\rr,\ss}; q,a),
\end{align}
where $\mathcal{H}_{[\lambda,\mu]}(T^{\rr,\ss}; q,a)$ is the 
composite HOMFLY-PT polynomial for $[\lambda,\mu]$ normalized 
by the condition $\mathcal{H}(U) = 1$ for the unknot $U$. 
\end{theorem}

{\it Proof.} This theorem formally results from the
coincidence of the $\tilde{J\!D}$\=polynomials in type $A$ with
the corresponding (reduced) Jones polynomials for torus knots
under the tilde-normalization. Generally, this claim is from 
Conjecture 2.1 in \cite{CJ}; it was verified there
for $A_{N-1}$ using the DAHA {\em shift operator\,}
(Proposition 2.3)
and papers \cite{LZ,Ste}. The weights were arbitrary there; we 
need them here for $[\la,\mu]_N$.\sq

\setcounter{equation}{0}
\section{\sc Examples and confirmations}\label{ExC}
We provide here examples of the composite 
DAHA-superpolynomials and discuss their symmetries. 
The first $5$ particular composite representations 
considered below are contained in the following table. 

\Yboxdim5pt 
$$
\begin{tabular}{ |c||c|c|c|c|c|}
\hline
$[b,c]$ & $[\omega_1,\omega_1]$ & $[\omega_1,\omega_2]$ 
& $[2\omega_1,\omega_1]$ & $\mathbf{[\omega_1,\omega_3]}$ 
& $[\omega_1+\omega_2,\omega_1]$\\ \hline
$[\lambda,\mu]$ & $[\yng(1),\yng(1)]$ & $[\yng(1),\yng(1,1)]$ 
& $[\yng(2),\yng(1)]$ & $[\yng(1),\yng(1,1,1)]$ 
& $[\yng(2,1),\yng(1)]$\\ \hline
$l$ & $2$ & $3$ & $2$ & $4$ & $3$\\ \hline\hline
$A_1$ & $\yng(2)$ & --- & $\yng(3)$ & --- & ---\\ \hline
$A_2$ & $\yng(2,1)$ & $\yng(2,2)$ & $\yng(3,2)$ & --- 
& $\yng(3,1)$\\ \hline
$A_3$ & $\yng(2,1,1)$ & $\yng(2,2,1)$ & $\yng(3,2,2)$ 
& $\yng(2,2,2)$ & $\yng(3,2,1)$\\ \hline
$A_4$ & $\yng(2,1,1,1)$ & $\yng(2,2,1,1)$ & $\yng(3,2,2,2)$ 
& $\yng(2,2,2,1)$ & $\yng(3,3,2,1)$\\ \hline
$A_5$ & $\yng(2,1,1,1,1)$ & $\yng(2,2,1,1,1)$ & $\yng(3,2,2,2,2)$ 
& $\yng(2,2,2,1,1)$ & $\yng(3,3,3,2,1)$\\ \hline
$A_6$ & $\yng(2,1,1,1,1,1)$ & $\yng(2,2,1,1,1,1)$ 
& $\yng(3,2,2,2,2,2)$ & $\yng(2,2,2,1,1,1)$ 
& $\yng(3,3,3,3,2,1)$\\ \hline
$\vdots$ & $\vdots$ & $\vdots$ & $\vdots$ & $\vdots$ 
& $\vdots$\\ \hline
\end{tabular}
$$

\subsection{\bf The adjoint representation}\label{adjoint}
The adjoint representation has the weight $\om_1+\om_n$
and is represented in our notation by the pair 
$[\omega_1,\omega_1] = [\yng(1),\yng(1)]$.  We consider 
this representation for two knots. 

\subsubsection{\sf Trefoil} 
 The adjoint DAHA 
superpolynomial for the trefoil is given by the formula 
\begin{multline*}
H\!D_{3,2}([\omega_1,\omega_1]; q,t,a) \;  = \;
1 + 2 qt + q^2 t^2 + 
 a(3 q^2 - q^3 + 2 qt^{-1} - q^2 t^{-1} - q^{3} t^{-1}\\ + 2 q^3 t) 
+ a^2(q^4 + q^2 t^{-2} - 2 q^3 t^{-2} + q^4t^{-2} + 2 q^3 t^{-1} - 
    2 q^4 t^{-1}) + 
 a^3(-q^4 t^{-3}\\ + q^5 t^{-3} + q^4 t^{-2} - q^5 t^{-2}).
\end{multline*}
Recall that it is defined by the relations
\begin{align}\label{specia}
H\!D_{\rr,\ss}([\la,\mu]; q,t,a\mapsto -t^{n+1}) 
= \widetilde{J\!D}^{A_n}_{\rr,\ss}(\la^\ast+\mu; q,t)
\end{align}
for $\la=\om_1,\mu=\om_1$ and all $n\geq 1$. 

The corresponding normalized adjoint HOMFLY-PT polynomial for 
the unframed trefoil is given by formula (2.17) from 
\cite{PBR}; see also Section \ref{Exs}. One has: 
\begin{multline*}
\mathcal{H}_{[\yng(1),\yng(1)]}(T^{3,2}) \; 
= \;  a^2(q^{-2} + q^2 + 2)
+ a^3(-2q^{-2} + q^{-1} + q - 2q^2 - 2)
 + a^4(q^{-2}\\ - 2q^{-1} - 2q + q^2 + 3)
+ a^5(q^{-1} + q - 2),
\end{multline*}
and we have the following confirmation of
Theorem \ref{CONNEC}:
$$a^2q^{-2}H\!D_{3,2}
([\omega_1, \omega_1]\,;\, q,t\mapsto q,a\mapsto -a) 
= \mathcal{H}_{[\yng(1),\yng(1)]}(T^{3,2}).$$
The super-duality from (\ref{comp-duality}) in this case is 
as follows: 
$$t^{-2}H\!D_{3,2}([\omega_1,\omega_1]; q,t,a) 
= q^{2}H\!D_{3,2}([\omega_1,\omega_1]; t^{-1}, q^{-1},a).$$
The evaluation formula (\ref{jones-sevalt}) reads
\begin{align*}
H\!D_{3,2}([\omega_1, \omega_1]; q,1,a) \;& = \;  
 \bigl(1 + q + aq\bigr)^2= \; H\!D_{3,2}(\omega_1; q,1,a)^2.
\end{align*}

\subsubsection{\sf The case of $T^{4,3}$}
The adjoint DAHA superpolynomial for the 
$(4,3)$\=torus knot $T(4,3)$ is given by the formula 
\begin{multline*}
H\!D_{4,3}([\omega_1,\omega_1]\,;\, q,t,a) 
\; = \;\\1+2 q t+2 q^2 t+3 q^2 t^2+2 q^3 t^2+q^4 t^2+4 q^3 t^3
+2 q^4 t^3+3 q^4 t^4+2 q^5 t^4+2 q^5 t^5+q^6 t^6\\+a (5 q^2
+5 q^3-q^4-3 q^5-2 q^6+2 qt^{-1}+q^2t^{-1}-q^3t^{-1}-q^4t^{-1}
-q^5t^{-1}+8 q^3 t\\+7 q^4 t+q^5 t-3 q^6 t-q^7 t+9 q^4 t^2
+7 q^5 t^2-q^6 t^2-q^7 t^2+8 q^5 t^3+5 q^6 t^3-q^7 t^3\\
+5 q^6 t^4+q^7 t^4+2 q^7 t^5)+a^2 (7 q^4+9 q^5-2 q^6-8 q^7
+q^2t^{-2}+2 q^3t^{-2}-2 q^4t^{-2}\\-3 q^5t^{-2}+q^7t^{-2}
+q^8t^{-2}+4 q^3t^{-1}+5 q^4t^{-1}-2 q^5t^{-1}-8 q^6t^{-1}
+q^8t^{-1}+8 q^5 t\\+9 q^6 t-2 q^7 t-3 q^8 t+7 q^6 t^2
+5 q^7 t^2-2 q^8 t^2+4 q^7 t^3+2 q^8 t^3+q^8 t^4)+a^3 (3 q^6
+5 q^7\\-q^8-3 q^9+q^4t^{-3}-q^5t^{-3}-3 q^6t^{-3}+3 q^7t^{-3}
+q^8t^{-3}-q^9t^{-3}+q^4t^{-2}+3 q^5t^{-2}\\-q^6t^{-2}
-8 q^7t^{-2}+4 q^8t^{-2}+q^9t^{-2}+2 q^5t^{-1}+5 q^6t^{-1}
-2 q^7t^{-1}-8 q^8t^{-1}+3 q^9t^{-1}\\+2 q^7 t+3 q^8 t-q^9 t
+q^8 t^2+q^9 t^2)+a^4 (q^9-q^7t^{-4}+2 q^8t^{-4}-q^9t^{-4}
+q^6t^{-3}\\-4 q^8t^{-3}+4 q^9t^{-3}-q^{10}t^{-3}+q^7t^{-2}
+q^8t^{-2}-4 q^9t^{-2}+2 q^{10}t^{-2}+q^8t^{-1}-q^{10}t^{-1})\\
+a^5 (-q^{10}t^{-5}+q^{11}t^{-5}-q^9t^{-4}+2 q^{10}t^{-4}
-q^{11}t^{-4}+q^9t^{-3}-q^{10}t^{-3}),
\end{multline*}
defined as for the trefoil.  Computed using (\ref{CompRJ}),
the corresponding normalized HOMFLY-PT polynomial is 
\begin{multline*}
\mathcal{H}_{[\yng(1),\yng(1)]}(T^{4,3}) \; = \;  
q^{-6}\bigl(
a^6(q^{12}+2 q^{10}+2 q^9+3 q^8+2 q^7+5 q^6+2 q^5+3 q^4
+2 q^3\\+2 q^2+1)
 + a^7(-2 q^{12}-q^{11}-4 q^{10}-4 q^9-6 q^8-4 q^7-8 q^6-4 q^5
-6 q^4-4 q^3
 -4 q^2\\-q-2) + a^8(q^{12}+2 q^{11}+2 q^{10}+2 q^9+5 q^8+2 q^7
+7 q^6+2 q^5
 +5 q^4+2 q^3+2 q^2+2 q\\+1) + a^9(-q^{11}-4 q^8+2 q^7-4 q^6+2 q^5
-4 q^4
 -q) + a^{10}(2 q^8-4 q^7+5 q^6-4 q^5\\+2 q^4)
+ a^{11}(2 q^7-4 q^6+2 q^5)\bigr).
\end{multline*}
We have the {\em connection formula\,} 
\begin{equation*}
a^6q^{-6}H\!D_{4,3}([\omega_1, \omega_1]\,;\, q,t\mapsto 
q,a\mapsto -a) = \mathcal{H}_{[\yng(1),\yng(1)]}(T^{4,3}).
\end{equation*}
The  {\em super-duality\,} reads  
$$t^{-6}H\!D_{4,3}([\omega_1, \omega_1]; q,t,a) 
= q^{6}H\!D_{4,3}([\omega_1, \omega_1]; t^{-1}, q^{-1},a),$$
and the evaluation at $t=1$ is as follows:
\begin{align*}
H\!D_{4,3}([\omega_1, \omega_1]; q,1,a) \;& = \; 
\bigl(1+q+2 q^2+q^3+a(q+2 q^2+2 q^3)+a^2 q^3\bigr)^2\\
& =  \; H\!D_{4,3}(\omega_1; q,1,a)^2.
\end{align*}

\subsection{\bf Column/row and a box}
Such diagrams correspond to the symmetric and wedge powers
of the fundamental representation.
\subsubsection{\sf Two-row and a box:}\,{$[2\omega_1,\omega_1] = 
[\yng(2),\yng(1)]$}.

Then the composite DAHA superpolynomial for the 
trefoil is 
\begin{multline*}
H\!D_{3,2}([2\omega_1,\omega_1]; q,t,a) \; = \;  1+q t+q^2 t
+q^3 t+q^3 t^2+2 q^4 t^2+q^5 t^3
+a (3 q^3\\+3 q^4-2 q^6-q^7+qt^{-1}+q^2t^{-1}
-q^4t^{-1}-q^5t^{-1}
+q^4 t+4 q^5 t+2 q^6 t-q^7 t+q^6 t^2\\+2 q^7 t^2)
+a^2 (2 q^6+4 q^7-q^8-2 q^9+q^3t^{-2}-q^5t^{-2}-q^6t^{-2}+q^8t^{-2}
+q^4t^{-1}\\+3 q^5t^{-1}+q^6t^{-1}-4 q^7t^{-1}-2 q^8t^{-1}
+q^9t^{-1}
+2 q^8 t+q^9 t)
+a^3 (q^{10}-q^7t^{-3}+q^9t^{-3}\\+q^6t^{-2}+q^7t^{-2}-2 q^8t^{-2}
-2 q^9t^{-2}
+q^{10}t^{-2}+q^{11}t^{-2}+2 q^8t^{-1}+q^9t^{-1}
-2 q^{10}t^{-1}\\-q^{11}t^{-1})
+a^4 (-q^{10}t^{-3}+q^{12}t^{-3}+q^{10}t^{-2}-q^{12}t^{-2}),
\end{multline*}
defined by  (\ref{specia}) for
$\la=2\omega_1, \mu=\omega_1$ and all $n\geq 1$.

The corresponding normalized HOMFLY-PT polynomial 
is given by formula (A.1) from \cite{PBR}, as well as computed 
using (\ref{CompRJ}).  It is
\begin{multline*}
\mathcal{H}_{[\yng(2),\yng(1)]}(T^{3,2}) \; = \; 
 q^{-3}\bigl(a^3(q^8 + 2q^6 + q^5 + q^4 + q^3 + q^2 + 1)
 + a^4(-1-q-2q^3\\-2q^4-q^5-2q^6-q^7-2q^9)
 + a^5(q+2q^4+q^5-q^6+2q^7+q^{10})
 + a^6(-q^5+q^6\\-2q^8+q^9) + a^7(-q^7+q^8+q^9-q^{10})\bigr),
\end{multline*}
and we have the relationship 
$$a^3q^{-3}H\!D_{3,2}([2\omega_1,\omega_1]\,;\, 
q,t\mapsto q,a\mapsto -a) = \mathcal{H}_{[\yng(2),\yng(1)]}
(T^{3,2}),$$
confirming Theorem \ref{CONNEC}.  The 
super-duality  here requires $[\omega_1,\omega_2]$,
which will be considered next. The evaluation at $t=1$  reads
\begin{align*}
&H\!D_{3,2}([2\omega_1,\omega_1]\,;\, q,1,a)  = 
\bigl(1\!+ q\! + aq\bigr)\times \bigl(1\!+q^2\!+q^3\!+q^4\!+a(q^2\\
&+q^3\!+q^4\!+q^5)+a^2 q^5\bigr)= H\!D_{3,2}(\omega_1; q,1,a)\times 
H\!D_{3,2}(2\omega_1; q,1,a).
\end{align*}

\subsubsection{\sf Two-column and a box:}\,
${[\omega_1,\omega_2] = 
[\yng(1),\yng(1,1)]}$. 

The DAHA superpolynomial for the trefoil reads 
\begin{multline*}
H\!D_{3,2}([\omega_1,\omega_2]\,;\, q,t,a) \; = \;  1+2 q t
+q t^2+q^2 t^2+q^2 t^3+q^2 t^4+q^3 t^5
+a (4 q^2\\-q^4+2 qt^{-2}-q^2t^{-2}-q^3t^{-2}+qt^{-1}
+2 q^2t^{-1}-2 q^3t^{-1}+q^2 t+3 q^3 t-q^4 t+3 q^3 t^2\\+q^4 t^3
+q^4 t^4)
+a^2 (3 q^4-q^5+q^2t^{-4}-2 q^3t^{-4}+q^4t^{-4}+2 q^2t^{-3}
-q^3t^{-3}-2 q^4t^{-3}\\+q^5t^{-3}+4 q^3t^{-2}-4 q^4t^{-2}
+2 q^3t^{-1}+q^4t^{-1}-q^5t^{-1}+q^4 t+q^5 t^2)
+a^3 (-q^4t^{-6}\\+q^5t^{-6}+q^3t^{-5}-2 q^4t^{-5}+q^5t^{-5}
+q^4t^{-4}-2 q^5t^{-4}+q^6t^{-4}+2 q^4t^{-3}-2 q^5t^{-3}\\+q^5t^{-2}
-q^6t^{-2}+q^5t^{-1})
+a^4 (-q^5t^{-7}+q^6t^{-7}+q^5t^{-5}-q^6t^{-5}),
\end{multline*}
where the specialization relations for all $n\geq 2$ are
$$H\!D_{3,2}([\omega_1,\omega_2]; q,t,a\mapsto -t^{n+1}) 
= \widetilde{J\!D}^{A_n}_{3,2}(\omega_2+\omega_n; q,t).$$
The corresponding normalized HOMFLY-PT polynomial 
is given by formula (A.4) from \cite{PBR}, as well as computed using 
(\ref{CompRJ}):
\begin{multline*}
\mathcal{H}_{ [\yng(1),\yng(1,1)]}(T^{3,2}) \; = \; 
 q^{-7}\bigl(a^3(q^2+2q^4+q^5+q^6+q^7+q^8+q^{10})
 + a^4(-2q-q^3-2q^4\\-q^5-2q^6-2q^7-q^9-q^{10})
 + a^5(1+2q^3-q^4+q^5+2q^6+q^9)
 + a^6(q-2q^2+q^4-q^5)\\ + a^7(-1+q+q^2-q^3)\bigr),
\end{multline*}
and we have the {\em connection formula\,} 
$$a^3q^{-5}H\!D_{3,2}([\omega_1,\omega_2]; 
q,t\mapsto q,a\mapsto -a) 
= \mathcal{H}_{ [\yng(1),\yng(1,1)]}(T^{3,2}).$$
The {\em super-duality\,} and {\em evaluation\,} are as follows:

\begin{align*} 
&t^{-3}H\!D_{4,3}([2\omega_1,\omega_1]; q,t,a) 
= q^{5}H\!D_{4,3}([\omega_1, \omega_2]; t^{-1}, q^{-1}, a),\\
&H\!D_{3,2}([\omega_1, \omega_2]; q,1,a) = 
\bigl(1 + q + aq\bigr)\bigl(1 + q + aq\bigr)^2\\
 &\ \ \ \ =  H\!D_{3,2}(\omega_1; q,1,a)\times 
H\!D_{3,2}(\omega_2; q,1,a).
\end{align*}

The corresponding standard superpolynomials
are   
\begin{align*}
&H\!D_{3,2}(\omega_1; q,t,a)=1+qt+aq,\\
H\!D_{3,2}(\om_2;q,t,a)&=
1\!+\frac{a^2 q^2}{t}\!+q t\!+q t^2\!+q^2 t^4 +
a \left(q\!+\frac{q}{t}\!+q^2 t\!+q^2 t^2\right).
\end{align*}
See e.g. \cite{CJ} and references therein.

\subsubsection{\sf Three-column and a box:}\,
${[\omega_1,\omega_3] = [\yng(1),\yng(1,1,1)]}$. 
This example is of deg$_{a}=5$, which matches our 
conjecture. The corresponding DAHA-superpolynomial 
for the trefoil is as follows: 

\begin{multline*}
H\!D_{3,2}([\omega_1,\omega_3]; q,t,a) \; = \;\\ 1+2 q t+q t^2
+q^2 t^2+q t^3+q^2 t^3+2 q^2 t^4+q^2 t^5+q^3 t^5+q^2 t^6+q^3 t^6
+q^3 t^7+q^3 t^9\\+q^4 t^{10}
+a (5 q^2+q^3-2 q^4+2 qt^{-3}-q^2t^{-3}-q^3t^{-3}+qt^{-2}
+2 q^2t^{-2}-2 q^3t^{-2}\\+qt^{-1}+3 q^2t^{-1}-q^3t^{-1}-q^4t^{-1}
+2 q^2 t+5 q^3 t-2 q^4 t+q^2 t^2+4 q^3 t^2+4 q^3 t^3+q^4 t^3
\\-q^5 t^3+q^3 t^4+3 q^4 t^4-q^5 t^4+q^3 t^5+2 q^4 t^5+3 q^4 t^6
+q^5 t^7+q^5 t^9)
+a^2 (q^3+6 q^4\\-3 q^5+q^2t^{-6}-2 q^3t^{-6}+q^4t^{-6}
+2 q^2t^{-5}-q^3t^{-5}-2 q^4t^{-5}+q^5t^{-5}+2 q^2t^{-4}\\
+q^3t^{-4}-4 q^4t^{-4}+q^5t^{-4}+q^2t^{-3}+5 q^3t^{-3}-5 q^4t^{-3}
+5 q^3t^{-2}-3 q^5t^{-2}+q^6t^{-2}\\+3 q^3t^{-1}+4 q^4t^{-1}
-4 q^5t^{-1}+4 q^4 t+q^5 t-q^6 t+2 q^4 t^2+2 q^5 t^2+3 q^5 t^3
-q^6 t^3\\+q^5 t^4+q^5 t^5+q^6 t^6)
+a^3 (q^5+q^6-q^7-q^4t^{-9}+q^5t^{-9}+q^3t^{-8}-2 q^4t^{-8}
+q^5t^{-8}\\+q^3t^{-7}-2 q^4t^{-7}+q^6t^{-7}+2 q^3t^{-6}-4 q^5t^{-6}
+2 q^6t^{-6}+4 q^4t^{-5}-5 q^5t^{-5}+q^6t^{-5}\\+3 q^4t^{-4}
-2 q^5t^{-4}-q^6t^{-4}+2 q^4t^{-3}+2 q^5t^{-3}-3 q^6t^{-3}+q^7t^{-3}
+4 q^5t^{-2}-2 q^6t^{-2}\\+2 q^5t^{-1}-q^6t^{-1}+q^6 t+q^6 t^2)
+a^4 (-q^5t^{-11}+q^6t^{-11}-q^5t^{-10}+q^6t^{-10}+q^4t^{-9}\\
-2 q^5t^{-9}+q^6t^{-9}+q^5t^{-8}-2 q^6t^{-8}+q^7t^{-8}+q^5t^{-7}
-2 q^6t^{-7}+q^7t^{-7}+2 q^5t^{-6}\\-2 q^6t^{-6}+q^6t^{-5}-q^7t^{-5}
+q^6t^{-4}-q^7t^{-4}+q^6t^{-3})
+a^5 (-q^6t^{-12}+q^7t^{-12}\\+q^6t^{-9}-q^7t^{-9}),
\end{multline*}
which is defined by (\ref{specia})
for all $n\geq 3$ and $\la=2\omega_1, \mu=\omega_3$: 
$$H\!D_{3,2}([\omega_1,\omega_3]; q,t,a\mapsto -t^{n+1}) 
= \widetilde{J\!D}^{A_n}_{3,2}(\omega_3+\omega_n; q,t).$$
The corresponding normalized HOMFLY-PT polynomial is 
\begin{multline*}
\mathcal{H}_{[\yng(1),\yng(1,1,1)]}(T^{3,2}) \; = \;  
q^{-16}\bigl(a^4(q^{20}+q^{18}+q^{16}+q^{15}+2 q^{14}+q^{13}
+2 q^{12}+q^{11}+2 q^{10}\\+q^9 +2 q^8+q^6) + 
a^5(-q^{20}-q^{18}-3 q^{16}-q^{15}-3 q^{14}-2 q^{13}
-4 q^{12}-2 q^{11} -4 q^{10}\\-2 q^9-4 q^8 -q^7-2 q^6-2 q^4) 
+ a^6(q^{18}+q^{16}+3 q^{14}+q^{13} +3 q^{12}+q^{11}
+3 q^{10}+2 q^9\\+3 q^8+q^7+ 2q^6+2 q^4+q^2) 
+ a^7(-q^3 - q^5 - q^7 - q^8 - q^9 - q^{10} 
- q^{12} - q^{14}) + a^8(q^7\\-q^6+q^5-q^4
 +2 q^3-q^2+q-1) + a^9(q^4-q^3-q+1)\bigr).
\end{multline*}
One has: $a^4q^{-10}H\!D_{3,2}([\omega_1,\omega_3]; 
q,t\mapsto q,a\mapsto -a) = 
\mathcal{H}_{[\yng(1,1,1),\yng(1)]}(T^{3,2})$ and
\begin{align*}
&H\!D_{3,2}([\omega_1,\omega_3]\,;\, q,1,a) =  
\bigl(1 + q + aq\bigr)\times \bigl(1 + q + aq\bigr)^3 \\
&\ \ \ \ \ \  \ =\ H\!D_{3,2}(\omega_1; q,1,a)\times 
H\!D_{3,2}(\omega_3; q,1,a).
\end{align*}

\subsection{\bf Three-hook and a box}\,
The last case is ${[\omega_1+\omega_2,\omega_1] = 
[\yng(2,1),\yng(1)]}$.

The corresponding DAHA-superpolynomial for the trefoil is 
\begin{multline*}
H\!D_{3,2}([\omega_1+\omega_2,\omega_1]\,;\, q,t,a) \; = \;\\  
1+3 q t-q t^2+4 q^2 t^2+q^3 t^2-2 q^2 t^3+4 q^3 t^3+q^4 t^3
-2 q^3 t^4
+4 q^4 t^4-q^4 t^5+3 q^5 t^5\\+q^6 t^6
+a (-2 q^2+12 q^3-4 q^5-q^6+2 qt^{-2}-q^2t^{-2}-q^3t^{-2}
+6 q^2t^{-1}
-2 q^3t^{-1}\\-2 q^4t^{-1}-q^5t^{-1}-4 q^3 t+16 q^4 t+q^5 t-4 q^6 t
-q^7 t-6 q^4 t^2+16 q^5 t^2-2 q^7 t^2\\-4 q^5 t^3
+12 q^6 t^3-2 q^7 t^3
-q^8 t^3-2 q^6 t^4+6 q^7 t^4-q^8 t^4+2 q^8 t^5)
+a^2 (-6 q^5+26 q^6\\-8 q^7-6 q^8+q^9+q^2t^{-4}-2 q^3t^{-4}
+q^4t^{-4}
+5 q^3t^{-3}-5 q^4t^{-3}-2 q^5t^{-3}+q^6t^{-3}\\
+q^7t^{-3}-q^3t^{-2}
+13 q^4t^{-2}-8 q^5t^{-2}-6 q^6t^{-2}+q^7t^{-2}
+q^8t^{-2}-4 q^4t^{-1}
+22 q^5t^{-1}\\-8 q^6t^{-1}-9 q^7t^{-1}+q^8t^{-1}+q^9t^{-1}-6 q^6 t
+22 q^7 t-8 q^8 t-2 q^9 t-4 q^7 t^2+13 q^8 t^2\\-5 q^9 t^2
+q^{10} t^2-q^8 t^3+5 q^9 t^3-2 q^{10} t^3+q^{10} t^4)
+a^3 (-3 q^8+14 q^9-9 q^{10}+q^{11}\\-q^4t^{-6}+q^5t^{-6}
+2 q^4t^{-5}-4 q^5t^{-5}+q^6t^{-5}+q^7t^{-5}+6 q^5t^{-4}-9 q^6t^{-4}
+3 q^8t^{-4}\\
-2 q^5t^{-3}+14 q^6t^{-3}-14 q^7t^{-3}-3 q^8t^{-3}+5 q^9t^{-3}
-3 q^6t^{-2}+21 q^7t^{-2}-18 q^8t^{-2}\\-3 q^9t^{-2}+3 q^{10}t^{-2}
-5 q^7t^{-1}+21 q^8t^{-1}-14 q^9t^{-1}+q^{11}t^{-1}-2 q^9 t
+6 q^{10} t\\-4 q^{11} t+q^{12} t+2 q^{11} t^2-q^{12} t^2)
+a^4 (q^{12}-q^{13}-q^6t^{-7}+q^7t^{-7}+q^6t^{-6}-3 q^7t^{-6}\\
+q^8t^{-6}+q^9t^{-6}+4 q^7t^{-5}-7 q^8t^{-5}+2 q^9t^{-5}
+2 q^{10}t^{-5}-q^{11}t^{-5}-q^7t^{-4}+8 q^8t^{-4}\\-11 q^9t^{-4}
+2 q^{10}t^{-4}+2 q^{11}t^{-4}-2 q^8t^{-3}+10 q^9t^{-3}
-11 q^{10}t^{-3}
+2 q^{11}t^{-3}\\+q^{12}t^{-3}-2 q^9t^{-2}+8 q^{10}t^{-2}
-7 q^{11}t^{-2}
+q^{12}t^{-2}-q^{10}t^{-1}+4 q^{11}t^{-1}-3 q^{12}t^{-1}\\
+q^{13}t^{-1})
+a^5 (-q^9t^{-7}+q^{10}t^{-7}+q^9t^{-6}-2 q^{10}t^{-6}+q^{11}t^{-6}
+2 q^{10}t^{-5}-3 q^{11}t^{-5}\\+q^{12}t^{-5}-q^{10}t^{-4}
+2q^{11}t^{-4}
-2 q^{12}t^{-4}+q^{13}t^{-4}+q^{12}t^{-3}-q^{13}t^{-3}),
\end{multline*}
defined by (\ref{specia}) for
$\la=\omega_1+\omega_2, \mu=\omega_1$ and all $n\geq 2$. 
The corresponding normalized HOMFLY-PT polynomial is 
\begin{multline*}
\mathcal{H}_{[\yng(2,1),\yng(1)]}(T^{3,2}) \; = \;  
q^{-8}\bigl(a^4(q^{14}+3 q^{12}-q^{11}+4 q^{10}-q^9+4 q^8-q^7
+4 q^6-q^5\\+3 q^4+q^2)
 + a^5(-2 q^{15}+q^{14}-5 q^{13}+4 q^{12}-10 q^{11}+5 q^{10}
-12 q^9+6 q^8-12 q^7\\
 +5 q^6-10 q^5+4 q^4 -5 q^3+q^2-2 q) + a^6(q^{16}-2 q^{15}
+6 q^{14} -6 q^{13}
+11 q^{12}-11 q^{11}\\+17 q^{10}-13 q^9+18 q^8-13 q^7+17 q^6 -11 
q^5+11 q^4 -6 q^3+6 q^2-2 q+1) + a^7(q^{16}\\-3 q^{15}+4 q^{14}-7 
q^{13}+10 q^{12}-14 q^{11} +14 q^{10}-18 q^9+18 q^8-18 q^7+14 q^6 
-14 q^5\\+10 q^4-7 q^3+4 q^2-3 q 
+1) + a^8(-q^{15}+2 q^{14}-3 q^{13}+5 q^{12}-7 q^{11}
+10 q^{10}-11 q^9\\+11 q^8
-11 q^7+10 q^6-7 q^5 +5 q^4-3 q^3+2 q^2-q) + a^9(q^{12}
-2 q^{11}+2 
q^{10} -3 q^9\\+4 q^8-3 q^7+2 q^6-2 q^5+q^4)\bigr), 
\end{multline*}
which reduces to the HOMFLY-PT polynomial as follows:
$$a^4q^{-6}H\!D_{3,2}([\omega_1+\omega_2,\omega_1]; 
q,t\mapsto q,a\mapsto -a) = 
\mathcal{H}_{[\yng(2,1),\yng(1)]}(T^{3,2}).$$  
The exact super-duality identity from (\ref{comp-duality}) is
$$t^{-6}H\!D_{3,2}([\omega_1+\omega_2,\omega_1]; q,t,a) 
= q^{6}H\!D_{3,2}([\omega_1+\omega_2,\omega_1]; t^{-1}, 
q^{-1},a).$$
The evaluation at $t=1$ from (\ref{jones-sevalt}) reads 
\begin{align*}
&H\!D_{3,2}([\omega_1+\omega_2,\omega_1]\,;\, q,1,a)\  = \  
\bigl(1 + q + aq\bigr)\\
&\times 
(1+q+a q) \left(1+q^2+q^3+q^4+a^2 q^5+
a \left(q^2+q^3+q^4+q^5\right)\right)\\
&=  H\!D_{3,2}(\omega_1; q,1,a)\times H\!D_{3,2}
(\omega_1\! +\! \omega_2; q,1,a).
\end{align*}

\subsection{\bf Two-rows and two-columns}
One of the two diagrams in the previous examples was
always a box. Let us discuss the cases when 
two-row and two-column diagrams are combined. 
They match well our conjectural formula (\ref{degaconj})
for deg${}_a$; we also checked directly the super-duality 
and other properties provided by the theorems above.

\subsubsection{\sf Two two-columns:\,} 
${[\omega_2,\omega_2] = 
[\yng(1,1),\yng(1,1)]}$.

$$H\!D_{3,2}([\om_2,\om_2]\,;\,q,t,a)=$$ 
\noindent
\smallskip
\renewcommand{\baselinestretch}{0.5} 
{\small
\(
1+a^6 \bigl(\frac{q^8}{t^{14}}-\frac{q^9}{t^{14}}-\frac{q^8}{t^{13}}
+\frac{q^{10}}{t^{13}}-\frac{q^8}{t^{12}}+\frac{2 q^9}{t^{12}}
-\frac{q^{10}}{t^{12}}+\frac{q^8}{t^{11}}-\frac{q^{10}}{t^{11}}
-\frac{q^9}{t^{10}}+\frac{q^{10}}{t^{10}}\bigr)+a^5 \bigl(
-\frac{q^6}{t^{13}}+\frac{q^7}{t^{13}}-\frac{q^6}{t^{12}}
+\frac{3 q^7}{t^{12}}-\frac{2 q^8}{t^{12}}+\frac{q^6}{t^{11}}
-\frac{q^7}{t^{11}}-\frac{q^8}{t^{11}}+\frac{q^9}{t^{11}}
+\frac{q^6}{t^{10}}-\frac{5 q^7}{t^{10}}+\frac{5 q^8}{t^{10}}
-\frac{q^9}{t^{10}}+\frac{2 q^8}{t^9}-\frac{2 q^9}{t^9}
+\frac{2 q^7}{t^8}-\frac{4 q^8}{t^8}+\frac{2 q^9}{t^8}
-\frac{q^8}{t^7}+\frac{q^9}{t^7}+\frac{q^8}{t^6}
-\frac{q^9}{t^6}\bigr)+a^4 \bigl(-\frac{2 q^5}{t^{11}}
+\frac{2 q^6}{t^{11}}+\frac{q^4}{t^{10}}-\frac{4 q^5}{t^{10}}
+\frac{4 q^6}{t^{10}}-\frac{q^7}{t^{10}}+\frac{q^6}{t^9}
-\frac{q^8}{t^9}+\frac{4 q^5}{t^8}-\frac{10 q^6}{t^8}
+\frac{7 q^7}{t^8}-\frac{q^8}{t^8}+\frac{2 q^5}{t^7}
-\frac{6 q^6}{t^7}+\frac{4 q^7}{t^7}+\frac{6 q^6}{t^6}
-\frac{9 q^7}{t^6}+\frac{3 q^8}{t^6}+\frac{3 q^6}{t^5}
-\frac{6 q^7}{t^5}+\frac{3 q^8}{t^5}+\frac{3 q^7}{t^4}
-\frac{3 q^8}{t^4}+\frac{2 q^7}{t^3}-\frac{2 q^8}{t^3}
+\frac{q^8}{t^2}\bigr)+2 q t+2 q t^2+q^2 t^2+2 q^2 t^3
+3 q^2 t^4+2 q^3 t^5+2 q^3 t^6+q^4 t^8+a^3 \bigl(q^6+q^7
-\frac{q^4}{t^9}+\frac{q^5}{t^9}+\frac{2 q^3}{t^8}
-\frac{5 q^4}{t^8}+\frac{3 q^5}{t^8}+\frac{2 q^3}{t^7}
-\frac{4 q^4}{t^7}+\frac{3 q^6}{t^7}-\frac{q^7}{t^7}
+\frac{6 q^4}{t^6}-\frac{12 q^5}{t^6}+\frac{7 q^6}{t^6}
-\frac{q^7}{t^6}+\frac{7 q^4}{t^5}-\frac{11 q^5}{t^5}
+\frac{2 q^6}{t^5}+\frac{2 q^7}{t^5}+\frac{q^4}{t^4}
+\frac{7 q^5}{t^4}-\frac{12 q^6}{t^4}+\frac{4 q^7}{t^4}
+\frac{10 q^5}{t^3}-\frac{12 q^6}{t^3}+\frac{2 q^7}{t^3}
+\frac{2 q^5}{t^2}+\frac{4 q^6}{t^2}-\frac{4 q^7}{t^2}
+\frac{7 q^6}{t}-\frac{5 q^7}{t}+2 q^7 t\bigr)
+a^2 \bigl(5 q^4+3 q^5-3 q^6+\frac{q^2}{t^6}-\frac{2 q^3}{t^6}
+\frac{q^4}{t^6}+\frac{4 q^2}{t^5}-\frac{4 q^3}{t^5}
-\frac{2 q^4}{t^5}+\frac{2 q^5}{t^5}+\frac{q^2}{t^4}
+\frac{4 q^3}{t^4}-\frac{9 q^4}{t^4}+\frac{3 q^5}{t^4}
+\frac{q^6}{t^4}+\frac{10 q^3}{t^3}-\frac{9 q^4}{t^3}
-\frac{2 q^5}{t^3}+\frac{q^6}{t^3}+\frac{4 q^3}{t^2}
+\frac{5 q^4}{t^2}-\frac{10 q^5}{t^2}+\frac{2 q^6}{t^2}
+\frac{15 q^4}{t}-\frac{10 q^5}{t}-\frac{q^6}{t}+10 q^5 t
-4 q^6 t+4 q^5 t^2-q^6 t^2+4 q^6 t^3+q^6 t^4\bigr)
+a \bigl(5 q^2+3 q^3-4 q^4+\frac{2 q}{t^3}-\frac{q^2}{t^3}
-\frac{q^3}{t^3}+\frac{2 q}{t^2}+\frac{q^2}{t^2}
-\frac{3 q^3}{t^2}+\frac{7 q^2}{t}-\frac{3 q^3}{t}
-\frac{2 q^4}{t}+10 q^3 t-3 q^4 t-q^5 t+6 q^3 t^2+q^4 t^2-q^5 t^2
+7 q^4 t^3-q^5 t^3+5 q^4 t^4-q^5 t^4+2 q^5 t^5+2 q^5 t^6\bigr).
\)
}
\renewcommand{\baselinestretch}{1.2} 
\smallskip

Note that the specializations $a=-t^{n+1}$ to $A_n$
begins here with $A_{n=3}$. We omit the formula
for $H\!D_{3,2}([2\om_1,2\om_1]\,;\,q,t,a)$, since it
can be readily obtained via the super-duality 
(checked numerically). Also,
$$H\!D_{3,2}([\om_2,\om_2]\,;\,q,t\!=\!1,a)=
(1+q +qa)^4.
$$

\subsubsection{\sf Two-column and two-row:\,} 
${[\omega_2, 2\omega_1] = 
[\yng(1,1),\yng(2)]}$.
Note that the $a$\=degree is $5$ in this example
vs. $6$ in the previous one; $\la\!\vee\!\mu$
contains now $3$ boxes (it is a $3$\=hook) in
(conjectural) formula (\ref{degaconj}). 
This formula is self-dual with respect
to $\,q\mapsto t^{-1},t\mapsto q^{-1}, a\mapsto a$
(up to $q^\bullet t^\bullet$).

$$H\!D_{3,2}([\om_2,2\om_1]\,;\,q,t,a)=$$ 
\noindent
\smallskip
\renewcommand{\baselinestretch}{0.5} 
{\small
\(
1+a^5 \bigl(-\frac{q^{11}}{t^7}+\frac{q^{13}}{t^7}
+\frac{q^{11}}{t^5}-\frac{q^{13}}{t^5}\bigr)
+a^4 \bigl(-\frac{q^8}{t^7}+\frac{q^{10}}{t^7}-\frac{q^9}{t^6}
-\frac{q^{10}}{t^6}+\frac{q^{11}}{t^6}+\frac{q^{12}}{t^6}
+\frac{q^7}{t^5}+\frac{q^8}{t^5}-\frac{q^9}{t^5}
-\frac{2 q^{10}}{t^5}+\frac{q^{12}}{t^5}+\frac{q^9}{t^4}
+\frac{q^{10}}{t^4}-\frac{2 q^{11}}{t^4}-\frac{q^{12}}{t^4}
+\frac{q^{13}}{t^4}+\frac{q^9}{t^3}+\frac{q^{10}}{t^3}
-\frac{q^{11}}{t^3}-\frac{q^{12}}{t^3}+\frac{q^{11}}{t^2}
-\frac{q^{13}}{t^2}+\frac{q^{11}}{t}\bigr)+q t+q^2 t+q^3 t
+q t^2+q^3 t^2+2 q^4 t^2+q^3 t^3+q^4 t^3+q^5 t^3+q^2 t^4
+q^5 t^4+q^4 t^5+q^5 t^5+q^6 t^6+a^3 \bigl(2 q^9+2 q^{10}
-q^{11}-\frac{q^7}{t^6}+\frac{q^9}{t^6}+\frac{q^4}{t^5}
-\frac{q^6}{t^5}-\frac{q^7}{t^5}-\frac{q^8}{t^5}+\frac{q^9}{t^5}
+\frac{q^{10}}{t^5}+\frac{2 q^6}{t^4}+\frac{q^7}{t^4}
-\frac{4 q^8}{t^4}-\frac{2 q^9}{t^4}+\frac{2 q^{10}}{t^4}
+\frac{q^{11}}{t^4}+\frac{q^5}{t^3}+\frac{2 q^6}{t^3}
+\frac{q^7}{t^3}-\frac{q^8}{t^3}-\frac{3 q^9}{t^3}
-\frac{2 q^{10}}{t^3}+\frac{q^{11}}{t^3}+\frac{q^{12}}{t^3}
+\frac{2 q^7}{t^2}+\frac{4 q^8}{t^2}-\frac{q^9}{t^2}
-\frac{4 q^{10}}{t^2}-\frac{q^{11}}{t^2}+\frac{q^7}{t}
+\frac{2 q^8}{t}+\frac{q^9}{t}+\frac{q^{10}}{t}-\frac{q^{11}}{t}
-\frac{q^{12}}{t}+q^9 t+q^{11} t^2\bigr)+a^2 \bigl(q^5+4 q^6
+4 q^7-q^9-q^{10}+\frac{q^3}{t^4}-\frac{q^5}{t^4}-\frac{q^6}{t^4}
+\frac{q^8}{t^4}+\frac{q^2}{t^3}+\frac{q^3}{t^3}-\frac{q^6}{t^3}
-\frac{3 q^7}{t^3}+\frac{2 q^9}{t^3}+\frac{3 q^4}{t^2}
+\frac{4 q^5}{t^2}-\frac{4 q^7}{t^2}-\frac{4 q^8}{t^2}
+\frac{q^{10}}{t^2}+\frac{2 q^4}{t}+\frac{q^5}{t}+\frac{4 q^6}{t}
+\frac{3 q^7}{t}-\frac{4 q^8}{t}-\frac{3 q^9}{t}+q^5 t+q^6 t+q^7 t
+4 q^8 t-q^{10} t+2 q^7 t^2+3 q^8 t^2+q^9 t^3+q^{10} t^3
+q^9 t^4\bigr)+a \bigl(q^2+3 q^3+2 q^4+2 q^5-q^6-2 q^7
+\frac{q}{t^2}+\frac{q^2}{t^2}-\frac{q^4}{t^2}-\frac{q^5}{t^2}
+\frac{q}{t}+\frac{2 q^3}{t}+\frac{2 q^4}{t}-\frac{q^5}{t}
-\frac{2 q^6}{t}-\frac{q^7}{t}+q^2 t+2 q^4 t+5 q^5 t
+2 q^6 t-q^7 t-q^8 t+3 q^4 t^2+2 q^5 t^2+2 q^6 t^2+2 q^7 t^2
-q^8 t^2+3 q^6 t^3+2 q^7 t^3
+q^5 t^4+q^6 t^4+q^8 t^4+q^7 t^5+q^8 t^5\bigr).
\)
}
\renewcommand{\baselinestretch}{1.2} 
\smallskip

The evaluation at $t=1$ from formula (\ref{jones-seval}) now
reads as follows:
\begin{align*}
&H\!D_{3,2}([\om_2,2\om_1]\,;\,q,t\!=\!1,a)\\ 
=
&(1 + q + a q)^2 (1 + q^2  + q^3 + q^4 + a (q^2+ q^3 + 
q^4 + q^5) + a^2 q^5,
\end{align*}
where the standard superpolynomial for $2\om_1$ is
$$
 H\!D_{3,2}(2\om_1\,;\,q,t,a)=
1\!+a^2 q^5\!+q^2 t\!+q^3 t\!+q^4 t^2+a 
\left(q^2\!+q^3\!+q^4 t\!+q^5 t\right).$$

Here and above we omit the 
formulas for the composite HOMFLY-PT polynomials; they 
do satisfy the Connection Theorem \ref{CONNEC}.

\setcounter{equation}{0}
\section{\sc Deligne-Gross series}\label{DG}
\subsection{\bf General procedure}
Here we consider the ``exceptional series": 
$$e\subset A_1 \subset A_2 \subset G_2 \subset D_4 
\subset F_4 \subset E_6 \subset E_7 \subset E_8,$$
discussed in \cite{DG}.  This is actually the bottom 
row of the triangle considered in that paper;
we are going to discuss it in full elsewhere.

Recall that the algebraic groups $G$ in 
this series are given a parameter $\nu$ in this paper
as follows: 
$$\nu(G) = \frac{h^\vee}{6},$$
where $h^{\vee}$ is the dual Coxeter number of $G$.
This very quantity provides the specializations
of our hyperpolynomials.
\smallskip

The {\em $E$\=hyperpolynomials\,} we will construct below
unify the DAHA-Jones polynomials (also called refined
polynomials) for $T^{3,2},T^{4,3}$
``colored" by the adjoint representation
for the groups of type $ADE$  in this series. The root
systems $G_2$ and $F_4$ play an important role in the
exceptional series, but we cannot incorporate them so far
(see also the end of this section).
\smallskip

As with the (colored) superpolynomial and 
hyperpolynomials of \cite{CJ,CJJ} and the present paper, 
this unification works by packaging the corresponding
DAHA-Jones polynomials into a single polynomial, denoted by
$H^{\mathfrak{ad}}_{\rr,\ss}(q,t,a)$, 
with an additional parameter $a$, where the individual 
polynomials are recovered via the following
specializations:
\begin{equation}\label{DGspec}
H^{\mathfrak{ad}}_{\rr,\ss}(q,t,a = -t^{\nu(G)}) = 
\widetilde{J\!D}^{G}_{\rr,\ss}(\mathfrak{ad}; q,t),
\hbox{\, excluding\, } G_2, F_4. 
\end{equation}
Thus $a$ is associated with the (dual) Coxeter number,  
rather than with the rank. Relations (\ref{DGspec})
appeared sufficient to determine 
$H^{\mathfrak{ad}}$ for $T^{3,2}$ and $T^{4,3}$, but 
this cannot be expected for arbitrary torus knots. 

In general, such polynomials cannot be uniquely determined via
these specializations for sufficiently complicated torus knots; 
one needs an infinite family of root systems in
(\ref{DGspec}) to restore $a$ for any knots.  
Practically speaking, however, only two specializations to $E_8$ and 
$E_7$ are enough for the trefoil.  We will demonstrate this in 
detail below.  Even more convincingly, the three specializations to 
$E_8$, $E_7$, and $E_6$ were enough for $T^{4,3}$;
the resulting polynomial has hundreds of terms.
\smallskip

Here we construct $H^{\mathfrak{ad}}_{\rr,\ss}$ for two knots, the 
trefoil $T^{3,2}$ and $T^{4,3}$. We will call this polynomial the 
\emph{adjoint exceptional hyperpolynomial}, since we consider only 
the adjoint representations.  As in \cite{CJJ}, we use the name 
``hyperpolynomial", since ``superpolynomial" is commonly
reserved for the root systems of type $A$.

For the trefoil we will show explicitly how 
$H^{\mathfrak{ad}}_{3,2}$ is obtained from the relevant DAHA-Jones 
polynomials for $E_8,E_7$ and the adjoint representation 
$\mathfrak{ad}$ whose highest weight is the highest short root 
$\vartheta$.  

For $T^{4,3}$, we obtain $H^{\mathfrak{ad}}_{4,3}$ using the same 
procedure, though $E_6$ is also required to find some
coefficients. Since the DAHA-Jones 
polynomials in these cases are rather long, we do not include them 
and instead refer the reader to \cite{CJJ} where they are posted.

Both $H^{\mathfrak{ad}}_{3,2}$ and $H^{\mathfrak{ad}}_{4,3}$ will 
satisfy all six of the defining specializations from (\ref{DGspec}), 
even though they are only constructed using two and three of these 
specializations, respectively. This is a convincing confirmation
that the formulas we found are meaningful. See Section  
\ref{DGspec2}, where we discuss this relations and some further
interesting symmetries.

\subsection{\bf E-type hyperpolynomials}
\subsubsection{\sf Trefoil}

Here we will demonstrate how $H^{\mathfrak{ad}}_{3,2}(q,t,a)$ is 
obtained from only the specializations (\ref{DGspec}) for 
$G$ of types $E_8,E_7$.  The relevant DAHA-Jones polynomial 
for $E_8$ from \cite{CJ} is
\begin{multline*}
\widetilde{J\!D}^{E_8}_{3,2}(\omega_8; q,t)  =  1 + q(t + t^6 
+t^{10} - t^{20} - t^{24} - t^{29})
 + q^2(t^{12} + t^{16} +t^{20}  - t^{26}\\ + t^{29} - 3t^{30} 
 - t^{34} - t^{35} - t^{39} + t^{44} + t^{49} + t^{53})
 + q^3(t^{29}+ t^{35} - t^{36} + t^{39} - t^{40} - t^{41}\\ 
 - t^{45} - 2t^{49} + t^{50} - t^{53} + t^{54}
 + t^{55} - t^{58} + 2t^{59} + t^{63} -t^{73})
 + q^4(t^{58} - t^{59} - t^{64} + t^{65}\\ - t^{68} + t^{69} 
 + t^{78} - t^{79} + t^{82} - t^{83})
+ q^5(-t^{87} + t^{88}),
\end{multline*}
and the relevant DAHA-Jones polynomial for $E_7$ is
\begin{multline*}
\widetilde{J\!D}^{E_7}_{3,2}(\omega_1; q,t)  =  1 + q(t + t^4 
+ t^6 - t^{12} - t^{14} - t^{17}) 
 + q^2(t^8 + t^{10} + t^{12}  -t^{16}\\ + t^{17} - 3t^{18} 
 - t^{20} - t^{21} - t^{23} +t^{26} + t^{29} + t^{31})
 + q^3(t^{17} + t^{21} - t^{22} + t^{23} - t^{24} - t^{25}\\ 
 - t^{27} - 2t^{29} +t^{30}  -t^{31} 
 +t^{32} +t^{33} - t^{34} +2t^{35} +t^{37} - t^{43})
 + q^4(t^{34} - t^{35} - t^{38} + t^{39}\\ - t^{40} + t^{41} + t^{46}
 - t^{47} + t^{48} - t^{49})
 + q^5(-t^{51} + t^{52}).
\end{multline*}
The (lexicographic) order in which these two polynomials are printed 
gives a perfect, one-to-one correspondence between their terms.  
Furthermore, this correspondence respects the signs $\pm$ of 
these terms.  

For example, in this correspondence $-q^2t^{39}$ in the 
$E_8$ polynomial 
is paired with $-q^2t^{23}$ in the $E_7$ polynomial. Determining the 
common exponent $x$ of $a$ that satisfies the right specializations 
from (\ref{DGspec}) readily reduces to finding a solution to
$39-5x=23-3x$, 
since $\nu(E_8)=5$ and $\nu(E_7)=3$.  Evidently, this solution 
is $x=8$, 
and the corresponding term in $H^{\mathfrak{ad}}_{3,2}$ will then 
be $-q^2t^{-1}a^8$.

Applying this procedure to every pair of terms in these two 
polynomials, 
the {\em adjoint exceptional hyperpolynomial\,} 
for the trefoil is 
\begin{multline*}
H^{\mathfrak{ad}}_{3,2}(q,t,a) \; = \;  1 + q(t - ta + a^2 - a^4 + 
t^{-1}a^5 - t^{-1}a^6) + q^2(t^2a^2 - ta^3 + a^4\\ + ta^5 + 
t^{-1}a^6 - 3a^6 + t^{-1}a^7 + a^7 - t^{-1}a^8 - t^{-1}a^9 + 
t^{-1}a^{10} - t^{-2}a^{11}) 
 + q^3(t^{-1}a^6\\ - a^7 + ta^7 + t^{-1}a^8 - a^8 - ta^8 
+ a^9 - 2t^{-1}a^{10} + a^{10}
+ t^{-2}a^{11} - t^{-1}a^{11} - a^{11} - t^{-2}a^{12}\\ + 
2t^{-1}a^{12} - t^{-2}a^{13} + t^{-2}a^{15}) 
 + q^4(t^{-2}a^{12} - t^{-1}a^{12} + t^{-1}a^{13} - a^{13} 
- t^{-2}a^{14}\\ + t^{-1}a^{14} + t^{-2}a^{16}
 - t^{-1}a^{16} - t^{-3}a^{17} + t^{-2}a^{17})
 + q^5(-t^{-3}a^{18} + t^{-2}a^{18}).
\end{multline*}

\subsubsection{\sf The case of $T^{4,3}$}
As it was mentioned above,
we will not provide the corresponding formulas for DAHA-Jones
polynomials for $E_{6,7,8}$ from \cite{CJ} here, since they 
are long.    
The adjoint exceptional hyperpolynomial for the torus knot $T^{4,3}$ 
can be constructed using essentially the same method as that for 
the trefoil. However, 
since the DAHA-Jones polynomials $\widetilde{J\!D}^{E_8}_{4,3}$ and 
$\widetilde{J\!D}^{E_7}_{4,3}$ have now different numbers of terms, 
their lexicographic orderings are (for some powers of $q$) 
insufficient to determine 
a correspondence between their respective monomials.  These 
few ambiguities are resolved by also considering 
$\widetilde{J\!D}^{E_6}_{4,3}$. 

Once such a correspondence between triples of monomials is 
established, the $a$\=degrees are uniquely restored using the 
relevant specializations from (\ref{DGspec}), as for the trefoil.  
The resulting hyperpolynomial is long, but we think that the
formula must be provided, since it has various symmetries
beyond those discussed in the paper and we expect that further
relations will be found. For instance, its connection to the root 
systems $F_4,G_2$ is an open problem. One has:
$H^{\mathfrak{ad}}_{4,3}(q,t,a)=$

\comment{
\begin{multline*}
H^{\mathfrak{ad}}_{4,3}(q,t,a) \; = \;\\  
1 + q(-t^{-1}a^6+t^{-1}a^5-a^4+a^2-t a+t)
+ q^2(-t^{-2}a^{11}+t^{-1}a^{10}-t^{-1}a^9\\
-t^{-1}a^8+t^{-1}a^7+a^7-4 a^6+t a^5+t^{-1}a^5+a^5-t a^4-t a^3
+t^2 a^2+t a^2+a^2-t^2 a\\-t a+t^2+t)
+ q^3(t^{-2}a^{15}-t^{-2}a^{13}+3t^{-1}a^{12}-3 t^{-1}a^{11}
-t^{-2}a^{11}-a^{11}+t^{-2}a^{10}\\+3 a^{10}-t^{-1}a^9-ta^8 
-t^{-1}a^8-3 a^8+4 t a^7+2 t^{-1}a^7+2 a^7-t^2 a^6-4 t a^6
+t^{-1}a^6\\-4 a^6+t^2 a^5+2 t a^5+a^5+a^4-t^3 a^3-t^2 a^3-2 t a^3
+t^3 a^2+2 t^2 a^2+t a^2-t^3 a-t^2 a\\+t^3) 
+ q^4(2 t^{-2}a^{17}-2 t^{-1}a^{16}-t^{-3}a^{16}+2 t^{-1}a^{15}
+t^{-2}a^{15}+2 t^{-1}a^{14}-t^{-2}a^{14}\\-2 t^{-1}a^{13}
-2 t^{-2}a^{13}-3 a^{13}+t a^{12}+3 t^{-1}a^{12}+6 a^{12}
-3 t a^{11}-5 t^{-1}a^{11}+t^{-2}a^{11}\\-4 a^{11}+2 t a^{10}
-2 t^{-1}a^{10}+2 a^{10}+t^2 a^9+2 t a^9+t^{-1}a^9+2 a^9
-4 t^2 a^8-5 t a^8\\+t^{-1}a^8-6 a^8+t^3 a^7+4 t^2 a^7+7 t a^7
+t^{-1}a^7-t^3 a^6-4 t^2 a^6-2 t a^6+t^{-1}a^6-t a^5\\+t^4 a^4
+t^3 a^4+3 t^2 a^4+t a^4+a^4-t^4 a^3-2 t^3 a^3-t^2 a^3-t a^3
+t^4 a^2+t^3 a^2+t^2 a^2)\\ 
+ q^5(-t^{-2}a^{21}+t^{-3}a^{20}+t^{-2}a^{19}-3 t^{-1}a^{18}
-t^{-3}a^{18}+3 t^{-1}a^{17}+4 t^{-2}a^{17}\\-t^{-3}a^{17}
+2 a^{17}-2 t^{-1}a^{16}-t^{-2}a^{16}-3 a^{16}+2 t^{-1}a^{15}
-2 t^{-2}a^{15}+3 t a^{14}+5 t^{-1}a^{14}\\-t^{-2}a^{14}+4 a^{14}
-t^2 a^{13}-6 t a^{13}-4 t^{-1}a^{13}-8 a^{13}+3 t^2 a^{12}
+6 t a^{12}-3 t^{-1}a^{12}\\+9 a^{12}-2 t^2 a^{11}-4 t a^{11}
-t^{-1}a^{11}+2 t^{-2}a^{11}+a^{11}-t^3 a^{10}-2 t^2 a^{10}
-4 t a^{10}\\-2 t^{-1}a^{10}-3 a^{10}+4 t^3 a^9+5 t^2 a^9
+7 t a^9+t^{-1}a^9+a^9-t^4 a^8-4 t^3 a^8-7 t^2 a^8\\-2 t a^8
+2 t^{-1}a^8-2 a^8+t^4 a^7+2 t^3 a^7+2 t^2 a^7+t a^7-2 a^7
+t^4 a^6+t^2 a^6+t a^6\\+t^{-1}a^6+a^6-t^5 a^5-t^4 a^5
-2 t^3 a^5-t^2 a^5-t a^5+t^5 a^4+t^4 a^4+t^3 a^4+t^2 a^4-t^4 a^3) 
\\+ q^6(-t^{-2}a^{23}+t^{-1}a^{22}+2 t^{-3}a^{22}-t^{-1}a^{21}
-2 t^{-2}a^{21}+t^{-3}a^{21}-t^{-1}a^{20}+t^{-2}a^{20}\\
+2 t^{-1}a^{19}+4 t^{-2}a^{19}-t^{-3}a^{19}+3 a^{19}-2 t a^{18}
-7 t^{-1}a^{18}+2 t^{-2}a^{18}-t^{-3}a^{18}-4 a^{18}\\+3 t a^{17}
+4 t^{-1}a^{17}-t^{-2}a^{17}-2 t^{-3}a^{17}+5 a^{17}
+4 t^{-1}a^{16}+t^{-3}a^{16}-2 a^{16}-3 t^2 a^{15}\\
-4 t a^{15}-t^{-1}a^{15}-3 t^{-2}a^{15}-6 a^{15}+t^3 a^{14}
+6 t^2 a^{14}+8 t a^{14}-2t^{-2} a^{14}+6 a^{14}\\
-3 t^3 a^{13}-4 t^2 a^{13}-9 t a^{13}+3 t^{-2}a^{13}+a^{13}
+2 t^2 a^{12}-t a^{12}-7 t^{-1}a^{12}+t^{-2}a^{12}\\-a^{12}
+t^4 a^{11}+3 t^3 a^{11}+2 t^2 a^{11}+3 t a^{11}+t^{-1}a^{11}
+t^{-2}a^{11}+3 a^{11}-4 t^4 a^{10}-2 t^3 a^{10}\\-6 t^2 a^{10}
-t a^{10}+t^{-1}a^{10}-2 a^{10}+t^5 a^9+t^4 a^9+4 t^3 a^9+2 t a^9
-3 a^9-t^3 a^8+2 t a^8\\+t^{-1}a^8+a^8-t^5 a^7-t^3 a^7-t a^7
-a^7+t^6 a^6+t^4 a^6+t^2 a^6+a^6) 
+ q^7(-t^{-3}a^{26}\\+t^{-1}a^{24}+t^{-3}a^{24}-t^{-1}a^{23}
-4 t^{-2}a^{23}+2 t^{-3}a^{23}-a^{23}+3 t^{-1}a^{22}
-t^{-2}a^{22}\\+t^{-3}a^{22}-t^{-4}a^{22}+a^{22}-t^{-1}a^{21}
+2 t^{-2}a^{21}+2 t^{-3}a^{21}+a^{21}-3 t a^{20}-5 t^{-1}a^{20}\\
+2 t^{-2}a^{20}-t^{-3}a^{20}-2 a^{20}+2 t^2 a^{19}+3 t a^{19}
+t^{-2}a^{19}-3 t^{-3}a^{19}+7 a^{19}-2 t^2 a^{18}\\-2 t a^{18}
+3 t^{-1}a^{18}+5 t^{-2}a^{18}-t^{-3}a^{18}-4 a^{18}-2 t^2 a^{17}
+2 t a^{17}-t^{-1}a^{17}-6 t^{-2}a^{17}\\+t^{-3}a^{17}-4 a^{17}
+3 t^3 a^{16}+2 t^2 a^{16}+5 t a^{16}+5 t^{-1}a^{16}
-3 t^{-2}a^{16}+a^{16}-t^4 a^{15}\\-3 t^3 a^{15}-3 t^2 a^{15}
-4 t a^{15}+2 t^{-1}a^{15}+t^{-2}a^{15}+t^4 a^{14}+5 t^2 a^{14}
-3 t a^{14}-7 t^{-1}a^{14}\\+t^{-2}a^{14}+t^4 a^{13}-t^3 a^{13}
+2 t^2 a^{13}-t a^{13}-t^{-1}a^{13}+t^{-2}a^{13}+7 a^{13}
-t^5 a^{12}-t^4 a^{12}\\-t^3 a^{12}-t^2 a^{12}-2 t a^{12}
-t^{-1}a^{12}+t^{-2}a^{12}-a^{12}+t^5 a^{11}+2 t^3 a^{11}
-t^2 a^{11}+t a^{11}\\-a^{11}-t^4 a^{10}+t^3 a^{10}-t^2 a^{10}
+2 t a^{10}+t^{-1}a^{10}-t^2 a^9-a^9+t a^8) + q^8(-t^{-3}a^{28}\\
+t^{-2}a^{27}-t^{-3}a^{27}+t^{-4}a^{26}-2 t^{-2}a^{25}
+t^{-3}a^{25}-a^{25}+t a^{24}+4 t^{-1}a^{24}-3 t^{-2}a^{24}\\
+t^{-3}a^{24}-t^{-4}a^{24}+t^{-1}a^{23}+5 t^{-3}a^{23}
-2 t^{-4}a^{23}-2 a^{23}-t a^{22}-2 t^{-1}a^{22}\\
-2 t^{-2}a^{22}+t^{-3}a^{22}+a^{22}+2 t^2 a^{21}-2 t^{-1}a^{21}
+3 t^{-2}a^{21}-3 t^{-3}a^{21}+4 a^{21}-t^3 a^{20}\\
-4 t a^{20}-t^{-1}a^{20}+6 t^{-2}a^{20}-2 t^{-3}a^{20}+2 a^{20}
+t^2 a^{19}+t a^{19}-5 t^{-1}a^{19}-t^{-2}a^{19}\\+t^{-3}a^{19}
-a^{19}+t^3 a^{18}-t^2 a^{18}+2 t a^{18}+6 t^{-1}a^{18}
-5 t^{-2}a^{18}-t^4 a^{17}-2 t^2 a^{17}+t a^{17}\\
+3 t^{-1}a^{17}-t^{-2}a^{17}+t^{-3}a^{17}-3 a^{17}
+t^3 a^{16}-t^{-1}a^{16}-2 a^{16}-t^3 a^{15}+t^2 a^{15}\\
-t^{-1}a^{15}+3 a^{15}-2 t a^{14}-t^{-1}a^{14}+t^{-2}a^{14}
+a^{14}-t^{-1}a^{13}+a^{13}+t^{-2}a^{12})
\\+ q^9(t^{-2}a^{29}-t^{-3}a^{29}-t^{-1}a^{28}+t^{-2}a^{28}
-t^{-3}a^{28}+t^{-4}a^{28}-2 t^{-3}a^{27}+2 t^{-4}a^{27}\\
+t^{-1}a^{26}-t^{-2}a^{26}+2 t^{-1}a^{25}-t^{-2}a^{25}
+3 t^{-3}a^{25}-2 t^{-4}a^{25}-2 a^{25}+t a^{24}+t^{-1}a^{24}\\
-5 t^{-2}a^{24}+5 t^{-3}a^{24}-t^{-4}a^{24}-a^{24}+2 t^{-1}a^{23}
-t^{-2}a^{23}-t^{-3}a^{23}-t a^{22}-t^{-1}a^{22}\\
+3 t^{-2}a^{22}-2 t^{-3}a^{22}+a^{22}+t^2 a^{21}-t a^{21}
-3 t^{-1}a^{21}+2 t^{-2}a^{21}+a^{21}-t^{-1}a^{20}\\
-t^{-2}a^{20}+2 a^{20}+t^{-1}a^{19}-a^{19}-t^{-2}a^{18}
+t^{-3}a^{18}) 
+ q^{10}(-t^{-3}a^{29}+2 t^{-4}a^{29}\\-t^{-5}a^{29}
+t^{-2}a^{28}-2 t^{-3}a^{28}+t^{-4}a^{28}-t^{-2}a^{26}
+2 t^{-3}a^{26}-t^{-4}a^{26}+t^{-1}a^{25}\\-2 t^{-2}a^{25}
+t^{-3}a^{25}) 
+ q^{11}(t^{-4}a^{30}-t^{-5}a^{30}).
\end{multline*}
}

\noindent
\smallskip
\renewcommand{\baselinestretch}{0.5} 
{\small
\(
1 + q\bigl(-t^{-1}a^6+t^{-1}a^5-a^4+a^2-t a+t\bigr)
+ q^2\bigl(-t^{-2}a^{11}+t^{-1}a^{10}-t^{-1}a^9
-t^{-1}a^8+t^{-1}a^7+a^7-4 a^6+t a^5+t^{-1}a^5+a^5-t a^4-t a^3
+t^2 a^2+t a^2+a^2-t^2 a-t a+t^2+t\bigr)
+ q^3\bigl(t^{-2}a^{15}-t^{-2}a^{13}+3t^{-1}a^{12}-3 t^{-1}a^{11}
-t^{-2}a^{11}-a^{11}+t^{-2}a^{10}+3 a^{10}-t^{-1}a^9-ta^8 
-t^{-1}a^8-3 a^8+4 t a^7+2 t^{-1}a^7+2 a^7-t^2 a^6-4 t a^6
+t^{-1}a^6-4 a^6+t^2 a^5+2 t a^5+a^5+a^4-t^3 a^3-t^2 a^3-2 t a^3
+t^3 a^2+2 t^2 a^2+t a^2-t^3 a-t^2 a+t^3\bigr) 
+ q^4\bigl(2 t^{-2}a^{17}-2 t^{-1}a^{16}-t^{-3}a^{16}+2 t^{-1}a^{15}
+t^{-2}a^{15}+2 t^{-1}a^{14}-t^{-2}a^{14}-2 t^{-1}a^{13}
-2 t^{-2}a^{13}-3 a^{13}+t a^{12}+3 t^{-1}a^{12}+6 a^{12}
-3 t a^{11}-5 t^{-1}a^{11}+t^{-2}a^{11}-4 a^{11}+2 t a^{10}
-2 t^{-1}a^{10}+2 a^{10}+t^2 a^9+2 t a^9+t^{-1}a^9+2 a^9
-4 t^2 a^8-5 t a^8+t^{-1}a^8-6 a^8+t^3 a^7+4 t^2 a^7+7 t a^7
+t^{-1}a^7-t^3 a^6-4 t^2 a^6-2 t a^6+t^{-1}a^6-t a^5+t^4 a^4
+t^3 a^4+3 t^2 a^4+t a^4+a^4-t^4 a^3-2 t^3 a^3-t^2 a^3-t a^3
+t^4 a^2+t^3 a^2+t^2 a^2\bigr) 
+ q^5\bigl(-t^{-2}a^{21}+t^{-3}a^{20}+t^{-2}a^{19}-3 t^{-1}a^{18}
-t^{-3}a^{18}+3 t^{-1}a^{17}+4 t^{-2}a^{17}-t^{-3}a^{17}
+2 a^{17}-2 t^{-1}a^{16}-t^{-2}a^{16}-3 a^{16}+2 t^{-1}a^{15}
-2 t^{-2}a^{15}+3 t a^{14}+5 t^{-1}a^{14}-t^{-2}a^{14}+4 a^{14}
-t^2 a^{13}-6 t a^{13}-4 t^{-1}a^{13}-8 a^{13}+3 t^2 a^{12}
+6 t a^{12}-3 t^{-1}a^{12}+9 a^{12}-2 t^2 a^{11}-4 t a^{11}
-t^{-1}a^{11}+2 t^{-2}a^{11}+a^{11}-t^3 a^{10}-2 t^2 a^{10}
-4 t a^{10}-2 t^{-1}a^{10}-3 a^{10}+4 t^3 a^9+5 t^2 a^9
+7 t a^9+t^{-1}a^9+a^9-t^4 a^8-4 t^3 a^8-7 t^2 a^8-2 t a^8
+2 t^{-1}a^8-2 a^8+t^4 a^7+2 t^3 a^7+2 t^2 a^7+t a^7-2 a^7
+t^4 a^6+t^2 a^6+t a^6+t^{-1}a^6+a^6-t^5 a^5-t^4 a^5
-2 t^3 a^5-t^2 a^5-t a^5+t^5 a^4+t^4 a^4+t^3 a^4+t^2 a^4
-t^4 a^3\bigr) 
+ q^6\bigl(-t^{-2}a^{23}+t^{-1}a^{22}+2 t^{-3}a^{22}-t^{-1}a^{21}
-2 t^{-2}a^{21}+t^{-3}a^{21}-t^{-1}a^{20}+t^{-2}a^{20}
+2 t^{-1}a^{19}+4 t^{-2}a^{19}-t^{-3}a^{19}+3 a^{19}-2 t a^{18}
-7 t^{-1}a^{18}+2 t^{-2}a^{18}-t^{-3}a^{18}-4 a^{18}+3 t a^{17}
+4 t^{-1}a^{17}-t^{-2}a^{17}-2 t^{-3}a^{17}+5 a^{17}
+4 t^{-1}a^{16}+t^{-3}a^{16}-2 a^{16}-3 t^2 a^{15}
-4 t a^{15}-t^{-1}a^{15}-3 t^{-2}a^{15}-6 a^{15}+t^3 a^{14}
+6 t^2 a^{14}+8 t a^{14}-2t^{-2} a^{14}+6 a^{14}
-3 t^3 a^{13}-4 t^2 a^{13}-9 t a^{13}+3 t^{-2}a^{13}+a^{13}
+2 t^2 a^{12}-t a^{12}-7 t^{-1}a^{12}+t^{-2}a^{12}-a^{12}
+t^4 a^{11}+3 t^3 a^{11}+2 t^2 a^{11}+3 t a^{11}+t^{-1}a^{11}
+t^{-2}a^{11}+3 a^{11}-4 t^4 a^{10}-2 t^3 a^{10}-6 t^2 a^{10}
-t a^{10}+t^{-1}a^{10}-2 a^{10}+t^5 a^9+t^4 a^9+4 t^3 a^9+2 t a^9
-3 a^9-t^3 a^8+2 t a^8+t^{-1}a^8+a^8-t^5 a^7-t^3 a^7-t a^7
-a^7+t^6 a^6+t^4 a^6+t^2 a^6+a^6\bigr) 
+ q^7\bigl(-t^{-3}a^{26}+t^{-1}a^{24}+t^{-3}a^{24}-t^{-1}a^{23}
-4 t^{-2}a^{23}+2 t^{-3}a^{23}-a^{23}+3 t^{-1}a^{22}
-t^{-2}a^{22}+t^{-3}a^{22}-t^{-4}a^{22}+a^{22}-t^{-1}a^{21}
+2 t^{-2}a^{21}+2 t^{-3}a^{21}+a^{21}-3 t a^{20}-5 t^{-1}a^{20}
+2 t^{-2}a^{20}-t^{-3}a^{20}-2 a^{20}+2 t^2 a^{19}+3 t a^{19}
+t^{-2}a^{19}-3 t^{-3}a^{19}+7 a^{19}-2 t^2 a^{18}-2 t a^{18}
+3 t^{-1}a^{18}+5 t^{-2}a^{18}-t^{-3}a^{18}-4 a^{18}-2 t^2 a^{17}
+2 t a^{17}-t^{-1}a^{17}-6 t^{-2}a^{17}+t^{-3}a^{17}-4 a^{17}
+3 t^3 a^{16}+2 t^2 a^{16}+5 t a^{16}+5 t^{-1}a^{16}
-3 t^{-2}a^{16}+a^{16}-t^4 a^{15}-3 t^3 a^{15}-3 t^2 a^{15}
-4 t a^{15}+2 t^{-1}a^{15}+t^{-2}a^{15}+t^4 a^{14}+5 t^2 a^{14}
-3 t a^{14}-7 t^{-1}a^{14}+t^{-2}a^{14}+t^4 a^{13}-t^3 a^{13}
+2 t^2 a^{13}-t a^{13}-t^{-1}a^{13}+t^{-2}a^{13}+7 a^{13}
-t^5 a^{12}-t^4 a^{12}-t^3 a^{12}-t^2 a^{12}-2 t a^{12}
-t^{-1}a^{12}+t^{-2}a^{12}-a^{12}+t^5 a^{11}+2 t^3 a^{11}
-t^2 a^{11}+t a^{11}-a^{11}-t^4 a^{10}+t^3 a^{10}-t^2 a^{10}
+2 t a^{10}+t^{-1}a^{10}-t^2 a^9-a^9+t a^8\bigr) + 
q^8\bigl(-t^{-3}a^{28}
+t^{-2}a^{27}-t^{-3}a^{27}+t^{-4}a^{26}-2 t^{-2}a^{25}
+t^{-3}a^{25}-a^{25}+t a^{24}+4 t^{-1}a^{24}-3 t^{-2}a^{24}
+t^{-3}a^{24}-t^{-4}a^{24}+t^{-1}a^{23}+5 t^{-3}a^{23}
-2 t^{-4}a^{23}-2 a^{23}-t a^{22}-2 t^{-1}a^{22}
-2 t^{-2}a^{22}+t^{-3}a^{22}+a^{22}+2 t^2 a^{21}-2 t^{-1}a^{21}
+3 t^{-2}a^{21}-3 t^{-3}a^{21}+4 a^{21}-t^3 a^{20}
-4 t a^{20}-t^{-1}a^{20}+6 t^{-2}a^{20}-2 t^{-3}a^{20}+2 a^{20}
+t^2 a^{19}+t a^{19}-5 t^{-1}a^{19}-t^{-2}a^{19}+t^{-3}a^{19}
-a^{19}+t^3 a^{18}-t^2 a^{18}+2 t a^{18}+6 t^{-1}a^{18}
-5 t^{-2}a^{18}-t^4 a^{17}-2 t^2 a^{17}+t a^{17}
+3 t^{-1}a^{17}-t^{-2}a^{17}+t^{-3}a^{17}-3 a^{17}
+t^3 a^{16}-t^{-1}a^{16}-2 a^{16}-t^3 a^{15}+t^2 a^{15}
-t^{-1}a^{15}+3 a^{15}-2 t a^{14}-t^{-1}a^{14}+t^{-2}a^{14}
+a^{14}-t^{-1}a^{13}+a^{13}+t^{-2}a^{12}\bigr)
+ q^9\bigl(t^{-2}a^{29}-t^{-3}a^{29}-t^{-1}a^{28}+t^{-2}a^{28}
-t^{-3}a^{28}+t^{-4}a^{28}-2 t^{-3}a^{27}+2 t^{-4}a^{27}
+t^{-1}a^{26}-t^{-2}a^{26}+2 t^{-1}a^{25}-t^{-2}a^{25}
+3 t^{-3}a^{25}-2 t^{-4}a^{25}-2 a^{25}+t a^{24}+t^{-1}a^{24}
-5 t^{-2}a^{24}+5 t^{-3}a^{24}-t^{-4}a^{24}-a^{24}+2 t^{-1}a^{23}
-t^{-2}a^{23}-t^{-3}a^{23}-t a^{22}-t^{-1}a^{22}
+3 t^{-2}a^{22}-2 t^{-3}a^{22}+a^{22}+t^2 a^{21}-t a^{21}
-3 t^{-1}a^{21}+2 t^{-2}a^{21}+a^{21}-t^{-1}a^{20}
-t^{-2}a^{20}+2 a^{20}+t^{-1}a^{19}-a^{19}-t^{-2}a^{18}
+t^{-3}a^{18}\bigr) 
+ q^{10}\bigl(-t^{-3}a^{29}+2 t^{-4}a^{29}-t^{-5}a^{29}
+t^{-2}a^{28}-2 t^{-3}a^{28}+t^{-4}a^{28}-t^{-2}a^{26}
+2 t^{-3}a^{26}-t^{-4}a^{26}+t^{-1}a^{25}-2 t^{-2}a^{25}
+t^{-3}a^{25}\bigr) 
+ q^{11}\bigl(t^{-4}a^{30}-t^{-5}a^{30}\bigr).
\)
}
\renewcommand{\baselinestretch}{1.2} 
\smallskip

\subsection{\bf Specializations}\label{DGspec2}
For $\{\rr,\ss\}\in\{\{3,2\},\{4,3\}\}$, the following 
specializations, which are special cases of (\ref{DGspec}),
are easily verified: 
$$H^{\mathfrak{ad}}_{\rr,\ss}(q,t,a = -t^{5}) = 
\widetilde{J\!D}^{E_8}_{\rr,\ss}(\omega_8; q,t),$$
$$H^{\mathfrak{ad}}_{\rr,\ss}(q,t,a = -t^{3}) = 
\widetilde{J\!D}^{E_7}_{\rr,\ss}(\omega_1; q,t),$$
$$H^{\mathfrak{ad}}_{\rr,\ss}(q,t,a = -t^{2}) = 
\widetilde{J\!D}^{E_6}_{\rr,\ss}(\omega_2; q,t),$$
$$H^{\mathfrak{ad}}_{\rr,\ss}(q,t,a = -t^{1}) = 
\widetilde{J\!D}^{D_4}_{\rr,\ss}(\omega_2; q,t),$$
$$H^{\mathfrak{ad}}_{\rr,\ss}(q,t,a = -t^{\frac{1}{2}}) 
= \widetilde{J\!D}^{A_2}_{\rr,\ss}(\omega_1 + \omega_2; q,t),$$
$$H^{\mathfrak{ad}}_{\rr,\ss}(q,t,a = -t^{\frac{1}{3}}) 
= \widetilde{J\!D}^{A_1}_{\rr,\ss}(2\omega_1; q,t).$$
The DAHA-Jones polynomials for the first four specializations 
may be found in \cite{CJJ}.  The last two DAHA-Jones polynomials
are specializations of the DAHA superpolynomials from Section 
\ref{adjoint}.

In addition to these defining specializations, the expressions for 
$H^{\mathfrak{ad}}_{\rr,\ss}$ possess two structures that resemble 
the ``canceling differentials" from \cite{DGR} and other papers.  
On the level of polynomials, these canceling differentials 
correspond to specializations of the parameters
with respect to which $H^{\mathfrak{ad}}_{\rr,\ss}$ 
becomes a single monomial.

The simplest such specialization corresponds to the evaluation 
at $t=1$ of DAHA-Jones polynomials.  On the level of 
hyperpolynomials, we set $a\mapsto -t^\nu= -1$, which readily
results in the relation
$$
H^{\mathfrak{ad}}_{\rr,\ss}(q,t=1,a=-1) = 1.
$$

The following example of a ``canceling differential"
is more interesting. We set $t=qa^6\,$.\,
Then 
$$H^{\mathfrak{ad}}_{3,2}(q,t,a) = q^3t^{-1}a^6 
+ (1-qt^{-1}a^6)\mathcal{Q}_{3,2}(q,t,a),$$
$$H^{\mathfrak{ad}}_{4,3}(q,t,a) = q^7t^{-1}a^6 
+ (1-qt^{-1}a^6)\mathcal{Q}_{4,3}(q,t,a)$$ 
for some 
polynomials $\mathcal{Q}_{\rr,\ss}(q,t,a)$.  Observe that 
$qt^{-1}a^6\mapsto -qt^{h^\vee - 1}$ in the specialization 
$a\mapsto -t^{\nu}$. Upon this specialization,
the above relations reflect the  
$PSL_{\,2}^\wedge(\Z)$\=invariance of the image of 
nonsymmetric Macdonald 
polynomials $E_{\vartheta}$ in the quotient of the polynomial 
representation of the corresponding DAHA 
under the relation $qt^{h^\vee - 1}=-1$ by its radical. 
However we did not check all details.
 
Let us also mention potential links of 
our hyperpolynomials evaluated at $a=-t^{-1}$ and $a=-1$ 
to the root systems $D_6$ and respectively 
$A_3$, which we are going to investigate elsewhere.


Finally, let us touch upon the root systems $G_2,F_4$
in the Deligne-Gross series.  For
$\nu(G_2) = \frac{2}{3}$ and for $\nu(F_4) = 
\frac{3}{2}$, the corresponding specializations of 
$H^{\mathfrak{ad}}_{\rr,\ss}$ resemble  the polynomials
$\widetilde{J\!D}^{G_2}_{3,2}(\omega_1; q,r,t)$ and 
$\widetilde{J\!D}^{F_4}_{3,2}(\omega_1; q,r,t)$
from \cite{CJ} at $r=t$, but do not coincide with them.
Hopefully, these specializations are connected with the
{\em untwisted\,} variants of these two DAHA-Jones
polynomials, but they are known so far only in the 
twisted setting.
\medskip

{\sf Conclusion.} Let us mention that we do not touch in this
paper the physics aspects of the composite superpolynomials
(and those for other root systems). See \cite{GJKS} concerning
the corresponding theory of {\em resolved conifold\,}; we 
thank Masoud Soroush for a discussion. In the 
{\em refined case\,} (related to open Gromov-Witten invariants), 
this approach reached so far only the simplest examples (our 
composite DAHA-superpolynomials are well ahead), but this is an 
important motivation of what we did in this paper. In contrast 
to conventional Gromov-Witten invariants, a systematic
theory of open Gromov-Witten invariants is not yet developed. 
See e.g. \cite{Ma} for a comprehensive account of this field. 

Finally, we note that the counterparts of the HOMFLY-PT 
polynomials for the classical series of root systems, for
instance Kauffman polynomials, can be generally addressed via 
Chern-Simons theory. Recall that DAHA provide a uniform 
theory of (refined) DAHA-Jones polynomials for any
root systems and arbitrary weights (for algebraic knots/links),
including the hyperpolynomials for the classical series 
(conjecturally for $B,C$). The exceptional DAHA-hyperpolynomials 
are quite a challenge for us (see the last section).


\comment{
To conclude, let us mention there exist some additional 
interesting specializations.  For example, 
$H^{\mathfrak{ad}}_{\rr,\ss}(q,t,a = 
-t^{-1})$ resembles 
$\widetilde{J\!D}^{D_6}_{\rr,\ss}(\omega_2; q,t)$, 
and $H^{\mathfrak{ad}}_{\rr,\ss}(q,t,a = -1)$ resembles 
$\widetilde{J\!D}^{A_3}_{\rr,\ss}(\omega_1 + \omega_3; q,t)$. 
}

\medskip
{\bf Acknowledgements.}
Our special thanks go to Sergei Gukov for his help and 
participation in this project, as well as for introducing 
the authors of this paper to each other. We thank Mikhail 
Khovanov, Aaron Lauda, Satoshi Nawata, Hoel Queffelec and 
David Rose for various clarifying discussions. We thank 
very much Masoud Soroush for his help with establishing
a connection with paper \cite{GJKS} (its part on the composite
Rosso-Jones formula), his providing details of the calculations
there and a general discussion. The first
author thanks Andras Szenes and the University of Geneva
for the invitation and hospitality.  The second author 
warmly thanks Sergei Gukov for his extensive patience
and advice over the past four years. We acknowledge
using the SAGE software  \cite{Sage-Combinat}
for the nonsymmetric Macdonald polynomials. 

\bibliographystyle{unsrt}

\end{document}